\documentclass[11pt]{article}
\usepackage[utf8]{inputenc}
\usepackage[english]{babel}
\usepackage{color,xcolor}
\usepackage{algpseudocode}
\usepackage{algorithm}
\usepackage{array}
\usepackage{mathtools}
\usepackage{amsmath}
\usepackage{amssymb}
\usepackage{graphicx}
\usepackage{caption}
\usepackage{subcaption}
\usepackage{amsfonts}
\usepackage{amsthm,amsmath,amssymb}

\usepackage{authblk}

\usepackage{geometry}
\usepackage{amsthm}
\newtheorem{theorem}{Theorem}[section]
\newtheorem{lemma}[theorem]{Lemma}

\newtheorem{proposition}[theorem]{Proposition}

\newtheorem{remark}[theorem]{Remark}

\newcommand{\dive}{\operatorname{div}}

\usepackage{datetime}

\date{\displaydate{date}}

\numberwithin{equation}{section}

\usepackage{scalerel,stackengine}
\stackMath
\newcommand\reallywidehat[1]{%
	\savestack{\tmpbox}{\stretchto{%
			\scaleto{%
				\scalerel*[\widthof{\ensuremath{#1}}]{\kern.1pt\mathchar"0362\kern.1pt}%
				{\rule{0ex}{\textheight}}%WIDTH-LIMITED CIRCUMFLEX
			}{\textheight}% 
		}{2.4ex}}%
	\stackon[-6.9pt]{#1}{\tmpbox}%
}

\begin{document}

\title{Doubling Inequality and Strong Unique Continuation for an Elliptic Transmission Problem \thanks{Finanziato dall’Unione europea - Next Generation EU, Missione 4 Componente 1 CUP  B53D23009200006}}
 \author[1]{Tianrui Dai}
 \author[1]{Elisa Francini}
 \author[1]{Sergio Vessella}

\affil[1]{Dipartimento di Matematica e Informatica, Università degli Studi di Firenze}

\date{}
\maketitle

\begin{abstract}

	We investigate the Strong Unique Continuation Property (SUCP) for elliptic equations with piecewise Lipschitz coefficients exhibiting jump discontinuities across a regular interface. We prove SUCP at the interface using a doubling inequality derived from a Carleman estimate with a singular weight. This result is intended as a first step toward solving the  inverse problem of  estimating the size of an unknown, merely measurable, inclusion inside a conductor from boundary measurements. 
	
		\medskip
	
	\noindent\textbf{Mathematics Subject Classification (2020)}
	Primary 35B60, 35J15, 35R30.
	
	\medskip
	
	\noindent \textbf{Keywords}
	 Unique Continuation Properties, Elliptic Equations, Inverse Problems.
\end{abstract}

\section{Introduction}
Investigating unique continuation properties and related quantitative estimates for Partial Differential Equations (PDEs) is central to inverse problems in PDEs, particularly with respect to uniqueness, data dependence, and size estimates of unknown inclusions. For the latter two aspects, strong unique continuation properties and their associated quantitative estimates have proven particularly effective. We therefore find it useful to briefly recall these topics here. 

Let $\Omega$ be an open connected set of $\mathbb{R}^n$, we say that the linear differential equation  

\begin{equation}\label{Intr-1}
	Lu =0 \quad \mbox{in } \Omega,
\end{equation}
 enjoys the Strong Unique Continuation Property (SUCP), if for any point $x_0\in \Omega$ and for any solution $u$ of \eqref{Intr-1} which satisfies

\begin{equation}\label{Intr-2}
	\int_{B_r(x_0)} u^2dx=\mathcal{O}\left(r^k\right), \ \ \  \mbox{as } r\rightarrow 0, \ \forall k\in \mathbb{N},
\end{equation}
then $u = 0$ in $\Omega$.

It is obvious that the SUCP implies the Weak Unique Continuation Property (WUCP),  that is, if  $\omega$ an open subset $\Omega$  then we have
$$ u =0 \mbox{ in } \omega\quad\Longrightarrow u = 0 \mbox{ in } \Omega.$$
For extended and detailed treatments of the aforementioned properties, we refer the reader to  \cite{HO63}, \cite{I98}, \cite{Ler}, \cite{Ves}.

Here we are interested in studying the SUCP (Strong Unique Continuation Property) for an equation of the form

\begin{equation*}
	Lu=\operatorname{div}\left(A(x)\nabla u\right)+b(x)\cdot \nabla u+c(x)u
\end{equation*}  
where $A$ is a symmetric, elliptic matrix-valued function in $\mathbb{R}^n$ (with real entries) that presents jump discontinuities across an interface $\Gamma$ (which is sufficiently regular) and is Lipschitz continuous on both sides of the interface $\Gamma$, moreover $b \in L^{\infty}(\mathbb{R}^n, \mathbb{R}^n)$, and $c \in L^{\infty}(\mathbb{R}^n)$. 

The interest in the issue of SUCP for the above operator arises, on the one hand, from the fact that, as shown in \cite{Mand}, \cite{Mil}, if $A$ is not Lipschitz continuous, even the WUCP may fail. On the other hand,  if $x_0 \in \Gamma$, whether the validity of \eqref{Intr-2} implies that $u$ vanishes is still an open problem \cite[Section 10.4.3]{Ler}, 
in spite of the fact that weaker forms of unique continuation across an interface have been established in \cite{DcFLVW}, \cite{FrLVW}, and \cite{Le-Lern}.

It is worth noting that this paper represents only an initial step in studying this topic, as we prove the following result: if $A(x) = a(x)I$ (with $I$ denoting the identity matrix on $\mathbb{R}^n$) and $a$ has jump discontinuities across a $C^{1,1}$ interface $\Gamma$, while being Lipschitz continuous on both sides of the interface, then, if $x_0 \in \Gamma$ and \eqref{Intr-2} hold, then $u$ vanishes identically. In fact, we prove this result as a consequence of a doubling inequality (Theorem \ref{Th-Doubling}), from which, as is well known \cite{Ga-Li}, the SUCP at $x_0$ follows. It is also important to recall that, by combining the doubling inequality with standard local regularity estimates for elliptic equations \cite{G-Tr}, it is easy to obtain — see again \cite{Ga-Li} — the property that $|u|^2$ is an $A_p$ weight (see \cite{Coi-Fe}). This latter property, in turn, has been insightfully reinterpreted and used in \cite{AMR2003}, \cite{ARS2000} to derive size estimates of unknown inclusions. Within this framework, we consider an inverse problem that we intend to address, for which the result presented in this paper represents a first step. 

Let $\Omega \subset \mathbb{R}^n$ be a conducting body (for $n=2,3$) with sufficiently regular boundary $\partial\Omega$, composed of two disjoint subdomains $\Omega_-$ and $\Omega_+$, separated by a known $C^{1,1}$ interface $\Gamma$. The conductivity is given by 
\begin{equation}\label{a}
a = a_-\chi_{\Omega_-} + a_+\chi_{\Omega_+},
\end{equation}
where $a_\pm(x)$ are Lipschitz continuous, so that $a$ has a jump discontinuity across $\Gamma$.

We assume the presence of an (unknown) inclusion $D \subset \Omega$, measurable in the Lebesgue sense, with modified conductivity $ka(x)$, where $k > 0$ and $k \neq 1$. Let $\phi \in H^{-1/2}(\partial\Omega)$ be a current density on the boundary, satisfying

$$
\int_{\partial\Omega} \phi \, dS = 0.
$$

The corresponding electrostatic potential $u$, defined up to a constant, solves the Neumann problem:

$$
\begin{cases}
	\operatorname{div}(a(1+(k-1)\chi_D)\nabla u) = 0, & \text{in } \Omega, \\
	a \dfrac{\partial u}{\partial \nu} = \phi, & \text{on } \partial \Omega, \\
	\int_{\Omega} u(x)\, dx = 0.
\end{cases}
$$

We aim to estimate the measure of $D$ from the boundary energy

$$
W = \int_{\partial \Omega} \phi u \, dS,
$$

by comparing it to

$$
W_0 = \int_{\partial \Omega} \phi u_0 \, dS,
$$
where $u_0$ solves the same problem in the absence of the inclusion $D$. This type of problem has been addressed in \cite{FrLVW} (and in \cite{Fr-Ve-Wa} for the case with complex coefficients); however, since the approach is based on three-sphere-type inequalities, estimating the measure of the inclusion $D$ involves imposing geometric assumptions that are stricter than simple Lebesgue measurability.

The proof of the  doubling inequality mentioned above begins with a suitable change of variables that flattens the interface (see \cite{A-E} for a similar strategy).
This transformation modifies the original equation
\begin{equation*}
%\label{eq1}
\dive \left(a \nabla u\right)=0  
\end{equation*}
with $a$ as in \eqref{a},  into the form
\begin{equation}
\label{eq2}
\dive \left(\tilde{a} g^{-1}\nabla u\right)=0  
\end{equation}
where $\tilde{a}$ is a real valued coefficient exhibiting a jump discontinuity across the interface $\{x_n=0\}$. The matrix-valued function $g$, arising from the flattening procedure,  satisfies the conditions \eqref{ellipticity}, \eqref{lipg},  \eqref{ginzero}, and \eqref{gatzero}, which will be detailed later.

For the transformed equation \eqref{eq2}, we derive a suitable Carleman estimate involving a singular weight centered at some point on the interface.   Carleman estimates -  introduced by Carleman in \cite{Car1} and \cite{Car2} - are a class  of weighted integral inequalities dependent  on a large parameter. For a comprehensive treatment of these estimates, see \cite{HO63}, \cite{Ler}, and \cite{Ves}. The specific technique used here, however, is inspired by the approach in \cite{E-V}.

The structure of the paper is as follows. In Section 2, we introduce the necessary notation and present a Carleman estimate for the operator  $\Delta_{g}=\dive\left(g^{-1}\nabla\right)$ multiplied by a function $\gamma$ that has a jump on the interface $\{x_n=0\}$. Section 3 is devoted to the proof of this estimate, which proceeds through several steps and departs from existing results in that it requires careful treatment of the interface behavior. Finally, in Section 4, we apply the Carleman estimate to establish a doubling inequality for solutions to equation \eqref{eq2}.

\section{Notation and main result for flat interface}
Let $n\geq2$ and $r>0$. We denote by  $B_{r}$
the open ball centered at $0$ with radius equals to $r$
in $\mathbb{R}^{n}$. We also write $B_{r}^{+}\coloneqq B_{r}\cap\left\{ x=\left(x^{\prime},x_{n}\right)|x^{\prime}\in\mathbb{R}^{n-1},x_{n}\geq0\right\} $,
and $B_{r}^{-}\coloneqq B_{r}\cap\left\{ x=\left(x^{\prime},x_{n}\right)|x^{\prime}\in\mathbb{R}^{n-1},x_{n}\leq0\right\} $.

Consider the second order operator 
$$\Delta_{g}\coloneqq\dive\left(g^{-1}\nabla\right),$$
where $g^{-1}=\left(g^{ij}\right)_{i,j=1}^{n}$ is
an elliptic symmetric matrix-valued function in $\mathbb{R}^{n}$,
that is, there exists some constant $\lambda\geq1$ such that for any
$\xi\in\mathbb{R}^{n}$, $x\in\mathbb{R}^{n}$, 
\begin{equation}\label{ellipticity}
\lambda^{-1}\left|\xi\right|^{2}\leq\left\langle g^{-1}\left(x\right)\xi,\xi\right\rangle \leq\lambda\left|\xi\right|^{2}.
\end{equation}
We denote by  $g=\left(g_{ij}\right)_{i,j=1}^{n}$ the inverse of $g^{-1}$
We  assume that $g^{-1}$ is Lipschitz continuous, that is there exist a positive $\Lambda$ such that for every $x,y\in\mathbb{R}^{n}$
 \begin{equation}\label{lipg}
\sum_{i,j=1}^{n}\left|g^{ij}(x)-g^{ij}(y)\right|\leq\Lambda|x-y|.\end{equation}

For any $\eta$ and $\xi\in\mathbb{R}^{n}$, we use the following
notation: the inner product 
\begin{equation*}\xi\cdot\eta\coloneqq\left\langle g\xi,\eta\right\rangle\end{equation*}
where $\left\langle ,\right\rangle $ is the general inner product
in $\mathbb{R}^{n}$,  
the norm
\begin{equation*}\left|\xi\right|^{2}\coloneqq\xi\cdot\xi=\left\langle g\xi,\xi\right\rangle\end{equation*} 
and the weighted gradient
\begin{equation*}\nabla_{g}\coloneqq g^{-1}\nabla\end{equation*}
Notice that
\begin{equation*}\Delta_{g}=\text{div}\left(\nabla_{g}\right).\end{equation*}

Due to the local character of the Carleman estimate, it is not restrictive to assume 
\begin{equation}\label{ginzero}
g^{-1}\left(0\right)=I_{n}=\left(\begin{array}{ccc}
1 & \cdots & 0\\
\vdots & \ddots & \vdots\\
0 & \cdots & 1
\end{array}\right).
\end{equation}

In order to manage the boundary terms that will appear on the interface, we need to make another assumption on our matrix $g$. This assumption, is not restricitve for the kind of application we have in mind, in which $g$ represents the change of variable needed to flatten the interface.
Let us assume that
\begin{equation}\label{gatzero}
 \left(g^{nk}\right)_{k=1}^{n-1}\left(x^{\prime},0\right)=0,\mbox{  and }
g^{nn}\left(x^{\prime},0\right)=1.
\end{equation}

Let us define the functional space 
\begin{equation*}
\mathcal{H}\left(B_{1}\right)\coloneqq\left\{ u\in H^1\left(B_{1} \right)\mbox{ such that }u|_{B_{1}^{+}}\in C^{\infty}\left(B_{1}^{+}\right)\text{ and }u|_{B_{1}^{-}}\in C^{\infty}\left(B_{1}^{-}\right)\right\}.
\end{equation*}
For any function $u\in \mathcal{H}\left(B_{1}\right)$  we use the following notation: $u^{\pm}=u|_{B_{1}^{\pm}}$.

%To simplify the notation, we will denote by  $\intop_{B_1}$
%to represent the sum of the integration in $B_{1}^{+}$ and %$B_{1}^{-}$, so that
%\[
%\intop_{B_1^\pm}u=\intop_{B_{1}^{+}}u^++\intop_{B_{1}^{-}}u^-.
%\]
%This will especially apply to integrals involving derivatives of $u$.
%
%We will only mantain the correct notation in the statements of our results.}

%We will also drop the underset $B_1^\pm$ for sake of readability,  when this will not cause confusion.

Let us denote by $\nu=-e_{n}$ and $\nu_{g}=g^{-1}\nu$. We will also denote by $\left[\cdot\right]$ 
the jump  on the inner boundary $\left\{ x_{n}=0\right\}$ that is
\begin{equation}\label{jump}
\left[z\right](x^\prime,0)=\lim_{x_n\to 0^+} z(x^\prime,x_n)-\lim_{x_n\to 0^-} z(x^\prime,x_n).
\end{equation}

Let us now introduce the following notation that will be useful to define the appropriate weight function in our Carleman estimate.
For some $\epsilon\in\left(0,1\right]$, let us set
\begin{equation}\label{psi}
\psi\left(s\right)\coloneqq s\operatorname{exp}\left(-\intop_{0}^{s}\frac{\text{dt}}{t^{1-\epsilon}\left(1+t^{\epsilon}\right)}\right),
\end{equation}
and
\begin{equation}\label{phi}
\phi\left(s\right)\coloneqq\frac{\psi\left(s\right)}{s\psi^{\prime}\left(s\right)}=1+s^{\epsilon}.
\end{equation}
We will denote by
\begin{equation}\label{sigma}
\sigma\left(x\right)\coloneqq\left|x\right|
\end{equation}
and
\begin{equation}\label{w}
w\coloneqq\psi\left(\sigma\right)
\end{equation}
will be called the weight function. 
Notice that 
\begin{equation}\label{asintpsi}
\frac{\psi(s)}{s}\to 1\mbox{ and } \psi'(s)\to 1 \mbox{ for } s\to 0^+
\end{equation}

Let us also set up the following notation:
for any  $v\in C^{2}\left(\mathbb{R}^{n}\setminus\left\{ 0\right\} \right)$
such that $v\left(x\right)>0$ and $\left|\nabla_{g}v\left(x\right)\right|>0$,  any  $f\in C^{\infty}\left(\mathbb{R}^{n}\setminus\left\{ 0\right\} \right)$ and  $\tau\geq1$, let us define

\begin{equation*}%\label{Pvtau}
P_{v,\tau}\left(f\right)\coloneqq v^{-\tau}\Delta_g\left(v^{\tau}f\right).\end{equation*}
We will use the following decomposition:
\begin{equation}\label{Pdecomposition}
P_{v,\tau}\left(f\right)=P_{v,\tau}^{\left(s\right)}\left(f\right)+2\tau \frac{\left|\nabla_{g}v\right|^{2}}{v^2}A_{v}\left(f\right) 
\end{equation}
where
\begin{equation}\label{Pvtaus}
P_{v,\tau}^{\left(s\right)}\left(f\right)\coloneqq\Delta_{g}f+\tau^{2}\frac{\left|\nabla_{g}v\right|^{2}}{v^{2}}f
\end{equation}
and
\begin{equation}\label{Av}
A_{v}\left(f\right)  =\frac{v}{\left|\nabla_{g}v\right|^{2}}\left(\nabla_{g}v\cdot\nabla_{g}f\right)+
\frac{1}{2}F_{v}^{g}f,
\end{equation}
with
\begin{equation}\label{Fgv}
F_{v}^{g}  =\frac{v\Delta_{g}v-\left|\nabla_{g}v\right|^{2}}{\left|\nabla_{g}v\right|^{2}} =\frac{v\Delta_{g}v}{\left|\nabla_{g}v\right|^{2}}-1.
\end{equation}
Notice that $P_{v,\tau}^{\left(s\right)}$ is the symmetric part, while $ \frac{\left|\nabla_{g}v\right|^{2}}{v^2}A_{v}\left(f\right) $ the anti-symmetric part, that is
\begin{equation*}%\label{anti}
\intop_{\mathbb{R}^{n}}\frac{\left|\nabla_{g}v\right|^{2}}{v^2}A_{v}\left(f\right)f=0 
\end{equation*}for every  $f\in C_0^{\infty}\left(\mathbb{R}^{n}\setminus\left\{ 0\right\} \right)$.

\begin{theorem}\label{teopiatta}
Let $g$ satisfy  assumptions \eqref{ellipticity}, \eqref{lipg},  \eqref{ginzero}, and \eqref{gatzero}, let $\gamma^+\in C^{0,1}(B_1^+)$ and $\gamma^-\in C^{0,1}(B_1^-)$ such that $\gamma^\pm\geq \gamma_0>0$ in $B_1^\pm$ and let $w$ as in \eqref{w}.
Then, there exist $\overline{r}\in(0,1)$ and $\tau_0>0$, depending only on $\lambda$, $\Lambda$,  $\gamma_0$,  $\|\gamma^\pm\|_{C^{0,1}(B_1^\pm)}$ and $\epsilon$ such that, for any $u\in \mathcal{H}\left(B_{1}\right)$ supported in $B^+_{\overline{r}}\setminus\{0\}$
and for any $\tau\geq \tau_0$ we have
\begin{align}\label{uconint}
&4\lambda\intop_{B_1}\gamma\sigma^2 w^{-2\tau}\left(\Delta_gu\right)^2 
%\geq\\
%&
\geq\frac{\tau\epsilon}{2(1+8\lambda^2)}\intop_{B_1}\gamma\sigma^\epsilon w^{-2\tau}\left|\nabla_{g}u\right|^{2}+
\frac{\tau^3\epsilon}{8\lambda}\intop_{B_1} \gamma \sigma^{\epsilon-2}w^{-2\tau}u^2+\nonumber
\\
&+\intop_{\{x_n=0\}}\!\!w^{-2\tau}\left[\left(4\tau w\left( \frac{\nabla_{g}w\cdot\nabla_{g}u}{|\nabla_{g}w|^2}\right) +\left(2\tau(n-2)\phi\left(\sigma\right)-4\tau^2\right)u\vphantom{\frac{l}{l}}\right)\left(\gamma\nabla_{g} u\cdot \nu_g\right)\right],%\nonumber
\end{align}
where $\gamma=\gamma^+\chi_{B_1^+}+\gamma^-\chi_{B_1^-}$ and $[\cdot]$ is defined as in \eqref{jump}. Notice that $\gamma\in C^{0,1}\left(B_1\setminus \{x_n=0\}\right)$, but is is not continuous on $\{x_n=0\}$.
\end{theorem}

\section{Proof of Theorem \ref{teopiatta}}
In this section, we prove the main result. First of all we find some estimates regarding the function $f=w^{-\tau}u$.
\subsection{Pointwise estimates for $f$}
This section contains the main calculations needed for the proof of our principal result.  
Consider a function $f\in \mathcal{H}\left(B_{1}\right)$.  Since we carry on the main calculations pointwise on each side of the interface, let us first restrict to $B_1^+$. In order to simplify the formulas, we will write $f$ even if the notation $f^+$ would be more correct.

A useful tool in the first part of the proof will be Rellich's identity. 
For any regular vector valued function $B$ and for any functions $\gamma$ and $f$, we have, almost everywhere in $B_1^+$, that
\begin{align}\label{eq:rellich_identity_result}
2\gamma\left(\Delta_{g}f\right)\left(B\cdot\nabla_{g}f\right)  &= 
2\dive\left(\left(\gamma B\cdot\nabla_{g}f\right)\nabla_{g}f\right)  -\dive\left(\gamma B\left|\nabla_{g}f\right|^{2}\right)+\nonumber \\
&+\gamma\left(\dive B\left|\nabla_{g}f\right|^{2}
-2\left(\partial_{j}B^{k}\right)g^{ji}\partial_{i}f\partial_{k}f
+B^{k}\partial_{k}g^{ij}\partial_{i}f\partial_{j}f\right)+\nonumber \\
&+\left(B\cdot\nabla_{g}\gamma\right)\left|\nabla_{g}f\right|^{2}  -2(\partial_{i}\gamma) B^{k}g^{ij}\partial_{i}f\partial_{k}f.
\end{align}

We will use the above identity with a special vector valued function given by
\begin{equation}\label{Bv}
B_v\coloneqq\frac{v\nabla_gv}{|\nabla_gv|^2},
\end{equation} 
for $v\in C^{2}\left(\mathbb{R}^{n}\setminus\left\{ 0\right\} \right)$
such that $v\left(x\right)>0$ and $\left|\nabla_{g}v\left(x\right)\right|>0$.
Let us now define
\begin{equation}\label{Mvg}
M_{v}^{g}\coloneqq S_{v}^{g}g
\end{equation}
where the symmetric matrix $S_{v}^{g}=\left\{ S_{v}^{g,ij}\right\} _{i,j=1}^{n}$ is given by
\begin{align}\label{Svg}
S_{v}^{g,ij}  \coloneqq&\frac{1}{2}\left(\left(\dive B_v-F_{v}^{g}\right)g^{ij}\right)-\left(\partial_{k}B_v^{j}\right)g^{ki} -\left(\partial_{k}B_v^{i}\right)g^{kj}+B_v^{k}\partial_{k}g^{ij}.
\end{align}
\begin{lemma}\label{lemmanuovo}
For any $v\in C^2(B_1^+\setminus\{0\})$ such that $v\left(x\right)>0$ and $\left|\nabla_{g}v\left(x\right)\right|>0$,  any $f\in  C^\infty\left(B_{1}^+\setminus\{0\}\right)$ and any $\gamma\in C^{0,1}\left(B_1^+\right)$, the following equality holds:
\begin{align}\label{eq:lemmanuovo}
\gamma\frac{v^{2}}{\left|\nabla_{g}v\right|^{2}}\left(P_{\tau}f\right)^{2}&=\, \gamma\frac{v^{2}}{\left|\nabla_{g}v\right|^{2}}\left(P_{\tau}^{\left(s\right)}f\right)^{2}+4\tau^{2}\gamma\frac{\left|\nabla_{g}v\right|^{2}}{v^{2}}\left(A_{v}\left(f\right)\right)^{2}+\\
&+ 4\tau\gamma M_{v}^{g}\nabla_{g}f\cdot\nabla_{g}f
 +2\tau(B_v\cdot\nabla_{g}\gamma)\left|\nabla_{g}f\right|^{2}+ \nonumber\\
 &-4\tau(\partial_{j}\gamma) B_v^{k}g^{ij}\partial_{i}f\partial_{k}f
 +\tau\gamma F_{v}^{g}\Delta_{g}\left(f^{2}\right)-2\tau^3\frac{\nabla_g v\cdot \nabla_g \gamma}{v}f^2+\nonumber\\
&+2\tau\dive\left(\tau^2 \frac{\gamma f^2\nabla_gv}{v}+2\left(\gamma B_v\cdot\nabla_{g}f\right)\nabla_{g}f  -\gamma \left|\nabla_{g}f\right|^{2}B_v\right)\nonumber
\end{align}
\end{lemma}
\begin{proof}
By \eqref{Pdecomposition}, we have
\begin{align}\label{dec}
\gamma\frac{v^{2}}{\left|\nabla_{g}v\right|^{2}}\left(P_{\tau}f\right)^{2}  &=\gamma\frac{v^{2}}{\left|\nabla_{g}v\right|^{2}}\left(P_{\tau}^{\left(s\right)}f\right)^{2}+4\tau^{2}\gamma\frac{\left|\nabla_{g}v\right|^{2}}{v^{2}}\left(A_{v}\left(f\right)\right)^{2}+\\
&+4\tau\gamma\left(\Delta_{g}f+\tau^{2}\frac{\left|\nabla_{g}v\right|^{2}}{v^{2}}f\right)A_{v}\left(f\right).\nonumber
\end{align}
Let us now deal with the last term in \eqref{dec} and let us write
\begin{equation}\label{MJ1J2}
4\tau\gamma\left(\Delta_{g}f+\tau^{2}\frac{\left|\nabla_{g}v\right|^{2}}{v^{2}}f\right)A_{v}\left(f\right)=J_1+J_2
\end{equation}
where
\[
J_{1}=4\tau\gamma\left(\Delta_{g}f\right)A_{v}\left(f\right),
\mbox{ and }
J_{2}=4\tau^{3}\gamma\frac{\left|\nabla_{g}v\right|^{2}}{v^{2}}fA_{v}\left(f\right).
\]
By \eqref{Av} and standard calculations
\begin{align}\label{calcJ2}
J_{2}  &=4\tau^{3}\gamma\frac{\left|\nabla_{g}v\right|^{2}}{v^{2}}f\left(
\frac{v}{\left|\nabla_{g}v\right|^{2}}\left(\nabla_{g}v\cdot\nabla_{g}f\right)+\frac{1}{2}F_{v}^{g}f\right)=\nonumber\\
&=
2\tau^{3}\dive\left(\frac{\gamma f^2\nabla_gv}{v}\right)
%
% \left[\gamma\frac{g^{jk}\partial_{k}v}{v}f^{2}\cdot\nu_{g}\right]
-2\tau^{3}f^{2}\dive\left(\frac{\gamma \nabla_gv}{v}\right)+\nonumber\\
 &\quad+2\tau^{3}\gamma\frac{v\left(\Delta_{g}v\right)-\left|\nabla_{g}v\right|^{2}}{v^{2}}f^{2}=
 \nonumber\\
 &=2\tau^{3}\dive\left(\frac{\gamma f^2\nabla_gv}{v}\right)-2\tau^3\frac{\nabla_g v\cdot \nabla_g \gamma}{v}f^2
\end{align}
Notice that we trivially have
\begin{equation}\label{eq:calculation_lapfsquare}
\Delta_{g}\left(f^{2}\right)
 =2\left|\nabla_{g}f\right|^{2}+2f\Delta_{g}f,
\end{equation}
and recalling again \eqref{Av} we get
\begin{align}\label{eq:calculation_I1}
J_{1} & =4\tau\gamma\left(\Delta_{g}f\right)\left(B_v\cdot\nabla_{g}f\right)-2\tau\gamma F_{v}^{g}\left|\nabla_{g}f\right|^{2}%\\
 %& 
 +\tau\gamma\Delta_{g}\left(f^{2}\right)F_{v}^{g}.%\nonumber 
\end{align}
By Rellich's identity \eqref{eq:rellich_identity_result} with $B=B_v$ we can write
\begin{align}\label{calcJ1}
J_{1} & =2\tau\dive\left(2\left(\gamma B_v\cdot\nabla_{g}f\right)\nabla_{g}f  -\gamma \left|\nabla_{g}f\right|^{2}B_v\right)+
 4\tau\gamma M_{v}^{g}\nabla_{g}f\cdot\nabla_{g}f\nonumber\\
 & +2\tau(B_v\cdot\nabla_{g}\gamma)\left|\nabla_{g}f\right|^{2} -4\tau(\partial_{j}\gamma) B_v^{k}g^{ij}\partial_{i}f\partial_{k}f%\\
 %& 
 +\tau\gamma F_{v}^{g}\Delta_{g}\left(f^{2}\right)
\end{align}
From \eqref{dec},  \eqref{MJ1J2}, \eqref{calcJ2} and \eqref{calcJ1} we finally get \eqref{eq:lemmanuovo}.
\end{proof}

\begin{lemma}
\label{prop:Lemma_of_calculation}
Let $v\in C^{2}\left(\mathbb{R}^{n}\setminus\left\{ 0\right\} \right)$
such that $v\left(x\right)>0$ and $\left|\nabla_{g}v\left(x\right)\right|>0$
for every $x\in\mathbb{R}^{n}\setminus\left\{ 0\right\} $. Let $S_{v}^{g},M_{v}^{g},F_{v}^{g}$ and $B_v$ be as in \eqref{Mvg},  \eqref{Svg} and \eqref{Bv}.  Let $\psi$, 
$\phi$ and $\sigma$ be as \eqref{psi}, \eqref{phi} and \eqref{sigma}.
We have 
\begin{equation}\label{le1}
M_{v}^{g}\nabla_{g}v=S_{v}^{g}\nabla v=0
\end{equation}

\begin{equation}\label{le2}
F_{\psi\left(v\right)}^{g}=\phi\left(v\right)F_{v}^{g}-\phi^{\prime}\left(v\right)v,
\end{equation}
and 
\begin{align}\label{le3}
M_{\psi\left(v\right)}^{g}\xi\cdot\eta & =v\phi^{\prime}\left(v\right)\left(\xi\cdot\eta-\frac{\left(\nabla_{g}v\cdot\xi\right)\left(\nabla_{g}v\cdot\eta\right)}{\left|\nabla_{g}v\right|^{2}}\right)+\phi\left(v\right)M_{v}^{g}\xi\cdot\eta,
\end{align}
for every $\xi,\eta\in\mathbb{R}^{n}$. Moreover,

\begin{equation}\label{le4}
F_{\sigma}^{g\left(0\right)}=n-2, \quad M_{\sigma}^{g\left(0\right)}=0
\end{equation}
\end{lemma}
\begin{proof}
These equalities follows by tedious but straigthforward calculations. See \cite{E-V} for a similar result.
\end{proof}

We will also need the following estimates regarding the weight function $w$, that follows immediately from  definitions \eqref{w} and \eqref{psi}:
\begin{lemma}
\label{prop:Lemma_intermedio}There exists a constant $C\geq 1$ depending only on $\lambda$, $\Lambda$ and $\epsilon$ such that
\begin{equation}\label{stimewl}
C^{-1}\sigma\leq w \leq \sigma,\quad C^{-1}\leq |\nabla_g w|\leq C
\end{equation}
\begin{equation*}\label{maggLap}
|\Delta_gw|\leq \frac{C}{w}
\end{equation*}

\begin{equation*}\label{Fgw}
|F_{w}^{g}|\leq C
\end{equation*}
More exactly,
\begin{equation}\label{gradienteprec}
\frac{\lambda^{-1}\sigma^{-2}}{4}\leq
\frac{\lambda^{-1}\sigma^{-2}}{(1+\sigma^\epsilon)^2}\leq \frac{|\nabla_g w|^2}{w^2}\leq \lambda \sigma^{-2}
\end{equation}
\end{lemma}

Let us now prove the following proposition, that will be an intermediate step toward our main Carleman estimate:

\begin{proposition}\label{propsenzaint}
Let $g$ satisfy assumptions \eqref{ellipticity}, \eqref{lipg} and \eqref{ginzero}, let $\gamma\in C^{0,1}(B_1^+)$ such that $\gamma\geq \gamma_0>0$ in $B_1^+$ and let $w$ as in \eqref{w}.
Then, there esist $\overline{r}\in(0,1)$ and $\tau_0>0$, depending only on $\lambda$, $\Lambda$,  $\gamma_0$,  $\|\gamma\|_{C^{0,1}(B_1^+)}$ and $\epsilon$ such that, for any $f\in C^\infty\left(B_1^+\setminus\{0\}\right)$ 
and for any $\tau\geq \tau_0$ we have
\begin{align}\label{senzaint}
\gamma\frac{w^{2}}{\left|\nabla_{g}w\right|^{2}}\left(P_{w,\tau}f\right)^{2}&\geq\frac{\gamma}{2}\frac{w^{2}}{\left|\nabla_{g}w\right|^{2}}\left(P_{w,\tau}^{\left(s\right)}f\right)^{2}+2\tau^{2}\gamma\frac{\left|\nabla_{g}w\right|^{2}}{w^{2}}\left(A_{w}\left(f\right)\right)^{2}+\nonumber\\
&+\frac{\tau\sigma^\epsilon\gamma\epsilon}{2}\left|\nabla_{g}f\right|^{2}+\tau^3\frac{ \gamma \lambda^{-1}\epsilon\sigma^{\epsilon-2}}{4}f^2+\dive G_2
\end{align}
almost everywhere in $B_{\overline{r}}^+$, where
\begin{align}\label{G2}
G_2=&2\tau^3 \frac{\gamma f^2\nabla_gw}{w}+4\tau\left(\gamma B_w\cdot\nabla_{g}f\right)\nabla_{g}f  -2\tau\gamma \left|\nabla_{g}f\right|^{2}B_w+\nonumber\\&+\tau(n-2)\gamma\phi\left(\sigma\right)\nabla_g(f^2).
\end{align}
\end{proposition}
\begin{proof} %{\color{magenta}Qui da qualche parte va collegato il supporto di $f$ con la piccolezza di $\sigma$}
We start from \eqref{eq:lemmanuovo} and estimate all the terms appearing there. First of all, let us denote by 
\begin{equation}\label{derivnormale}\nabla_{g}^{N}f\coloneqq\left(\nabla_{g}\sigma\cdot\nabla_{g}f\right)\frac{\nabla_{g}\sigma}{\left|\nabla_{g}\sigma\right|^{2}}=
\left(\nabla_{g}w\cdot\nabla_{g}f\right)\frac{\nabla_{g}w}{\left|\nabla_{g}w\right|^{2}},
\end{equation}
which is the projection of $\nabla_{g}f$ on the normal direction $\nabla_{g}\sigma$.

Let us also denote by 
\begin{equation*}%\label{derivtang}
\nabla_{g}^{T}f\coloneqq\nabla_{g}f-\nabla_{g}^{N}f
\end{equation*}
the tangential part.

We trivially have 
\begin{equation}\label{modulgrad}
\left|\nabla_{g}^{T}f\right|^{2}+\left|\nabla_{g}^{N}f\right|^{2}=\left|\nabla_{g}f\right|^{2}.
\end{equation}
Thanks to \eqref{le3} (with $v=\sigma$ and $\xi=\eta=\nabla_g f$) and by \eqref{modulgrad}  we have
\begin{equation}\label{inter1}
M_{w}^{g}\nabla_{g}f\cdot\nabla_{g}f  =\sigma\phi^{\prime}\left(\sigma\right)\left|\nabla_{g}^{T}f\right|^{2}
 +\phi\left(\sigma\right)M_{\sigma}^{g}\nabla_{g}f\cdot\nabla_{g}f
\end{equation} 
 Notice that, by using \eqref{le1} and \eqref{derivnormale} with $v=\sigma$ we get
 \begin{equation}\label{inter2}
M_{\sigma}^{g}\nabla_{g}^{N}f  =\frac{\left(\nabla_{g}\sigma\cdot\nabla_{g}f\right)}{\left|\nabla_{g}\sigma\right|^{2}}M_{\sigma}^{g}\nabla_{g}\sigma=0
\end{equation} 
Hence, from \eqref{inter1},  \eqref{inter2} and the second equality in \eqref{le4}, we have
 \begin{equation}\label{inter3}
M_{w}^{g}\nabla_{g}f\cdot\nabla_{g}f  =\sigma\phi^{\prime}\left(\sigma\right)\left|\nabla_{g}^{T}f\right|^{2}+\phi\left(\sigma\right)\left(M_{\sigma}^{g}-M_{\sigma}^{g\left(0\right)}\right)\nabla_{g}^{T}f\cdot\nabla_{g}^{T}f
\end{equation} 
On the other hand, a direct calculation and assumptions \eqref{ellipticity} and \eqref{lipg} on $g$ give that there is a constant $C_0$ depending only on $\lambda$ and $\Lambda$ such that 
 \begin{equation}\label{lipM}
\left|M_{\sigma}^{g}-M_{\sigma}^{g\left(0\right)}\right|\leq C_0\sigma
\end{equation}
and, hence, by \eqref{inter3} and \eqref{lipM} we finally have
\begin{equation}\label{inter4}
M_{w}^{g}\nabla_{g}f\cdot\nabla_{g}f  \geq\sigma\left(\phi^{\prime}\left(\sigma\right)-C_0\phi\left(\sigma\right)\right)\left|\nabla_{g}^{T}f\right|^{2}.
\end{equation}
Let us notice that, since $\epsilon\in(0,1)$,  if we take $\delta\in(0,1)$ that will be chosen later on, we have
\begin{equation*}%\label{inter4.5}
\sigma\left(\phi^{\prime}\left(\sigma\right)-C_0\phi\left(\sigma\right)\right)=\sigma\left(\epsilon \sigma^{\epsilon-1}-C_0(1+\sigma^\epsilon)\right)\geq (1-\delta)\epsilon\sigma^\epsilon
\end{equation*}
if  
\begin{equation}\label{inter4.6}
\sigma\leq \overline{r}_0\coloneqq\min\left\{1,\left(\frac{\delta\epsilon}{2C_0}\right)^{\frac{1}{1-\epsilon}}\right\}.
\end{equation}
Notice that $\overline{r}_0$ depends only on $\delta$, $\Lambda$,  $\lambda$ and $\epsilon$.

By \eqref{inter4} we can write
\begin{equation}\label{inter5}
M_{w}^{g}\nabla_{g}f\cdot\nabla_{g}f  \geq \frac{(1-\delta)\epsilon\sigma^\epsilon}{1}\left|\nabla_{g}^Tf\right|^{2}= \frac{(1-\delta)\epsilon\sigma^\epsilon}{1}\left(\left|\nabla_{g}f\right|^{2}-\left|\nabla_{g}^{N}f\right|^{2}\right),
\end{equation}
and, by \eqref{derivnormale} and \eqref{Av}, we have
\begin{equation}\label{inter6}
\left|\nabla_{g}^{N}f\right|^{2}=\frac{|\nabla_gw|^2}{w^2}\left(A_w(f)-\frac{1}{2}F_w^gf\right)^2\leq 2\frac{|\nabla_gw|^2}{w^2}\left((A_w(f))^2+\frac{1}{4}(F_w^gf)^2\right).
\end{equation}
From \eqref{inter5}, \eqref{inter6} and  \eqref{Fgw}, we get
%\begin{equation}\label{inter7}
%M_{w}^{g}\nabla_{g}f\cdot\nabla_{g}f  \geq \frac{(1-\delta)\epsilon\sigma^\epsilon}{1}\left|\nabla_{g}f\right|^{2}- 2(1-\delta)\epsilon\sigma^\epsilon\frac{|\nabla_gw|^2}{w^2}\left((A_w(f))^2+\frac{1}{4}(F_w^gf)^2\right)
%\end{equation}
%By \eqref{Fgw},
\begin{equation}\label{inter8}
M_{w}^{g}\nabla_{g}f\cdot\nabla_{g}f  \geq \frac{(1-\delta)\epsilon\sigma^\epsilon}{1}\left|\nabla_{g}f\right|^{2}- 2(1-\delta)\epsilon\sigma^\epsilon\frac{|\nabla_gw|^2}{w^2}\left((A_w(f))^2+Cf^2\right).
\end{equation}
Recalling  \eqref{asintpsi} we also have that 
\begin{equation}\label{oooo}
 \left|2\tau(B_w\cdot\nabla_{g}\gamma)\left|\nabla_{g}f\right|^{2}
-4\tau(\partial_{j}\gamma) B_w^{k}g^{ij}\partial_{i}f\partial_{k}f\right|\leq C \tau \sigma\left|\nabla_{g}f\right|^{2}
\end{equation}
where $C$ depends on  $\lambda$ and $\left\Vert \gamma\right\Vert _{C^{0,1}}$.

Now,  by using \eqref{eq:lemmanuovo}, (with $v=w$), \eqref{inter8} and \eqref{oooo} and by setting
\begin{equation}\label{G1}
G_1=\tau^2 \frac{\gamma f^2\nabla_gw}{w}+2\left(\gamma B_w\cdot\nabla_{g}f\right)\nabla_{g}f  -\gamma \left|\nabla_{g}f\right|^{2}B_w
\end{equation}
we have, almost everywhere in $B_1^+$, 
\begin{align}\label{inter9}
\gamma&\frac{w^{2}}{\left|\nabla_{g}w\right|^{2}}\left(P_{w,\tau}f\right)^{2}\geq\, \gamma\frac{w^{2}}{\left|\nabla_{g}w\right|^{2}}\left(P_{w,\tau}^{\left(s\right)}f\right)^{2}+4\tau^{2}\gamma\frac{\left|\nabla_{g}w\right|^{2}}{w^{2}}\left(A_{w}\left(f\right)\right)^{2}+\nonumber\\
&+ 4\tau\gamma \left( (1-\delta)\epsilon\sigma^\epsilon\left|\nabla_{g}f\right|^{2}-2(1-\delta) \epsilon\sigma^\epsilon\frac{|\nabla_gw|^2}{w^2}\left((A_w(f))^2+Cf^2\right)\right)+
 \nonumber\\
 &-C \tau \sigma\left|\nabla_{g}f\right|^{2}
 +\tau\gamma F_{w}^{g}\Delta_{g}\left(f^{2}\right)-2\tau^3\frac{\nabla_g w\cdot \nabla_g \gamma}{w}f^2+2\tau\dive G_1=\nonumber\\
&=\gamma\frac{w^{2}}{\left|\nabla_{g}w\right|^{2}}\left(P_{w,\tau}^{\left(s\right)}f\right)^{2}+
\tau\sigma^\epsilon\left(4(1-\delta)\gamma\epsilon-C(1-\delta)\sigma^{1-\epsilon}\right)\left|\nabla_{g}f\right|^{2}+\nonumber\\
&+4\tau^{2}\gamma\frac{\left|\nabla_{g}w\right|^{2}}{w^{2}}\left(1-\frac{2(1-\delta)\epsilon\sigma^\epsilon}{\tau}\right)\left(A_{w}\left(f\right)\right)^{2}+\tau\gamma F_{w}^{g}\Delta_{g}\left(f^{2}\right)+\nonumber\\
&-8(1-\delta)\tau\gamma\epsilon\sigma^\epsilon\frac{|\nabla_gw|^2}{w^2}Cf^2
 -2\tau^3\frac{\nabla_g w\cdot \nabla_g \gamma}{w}f^2+2\tau\dive G_1
\end{align}
Now let us consider the term $ \tau\gamma F_{w}^{g}\Delta_{g}\left(f^{2}\right)$.
From \eqref{le2} and \eqref{le4},  we have
\begin{equation}\label{inter10}
F_{w}^{g}=\left(n-2\right)\phi\left(\sigma\right)-H\left(\sigma\right),
\end{equation}
where
\begin{equation}\label{H}
H\left(\sigma\right)=\sigma\phi^{\prime}\left(\sigma\right)-\phi\left(\sigma\right)\left(F_{\sigma}^{g}-F_{\sigma}^{g\left(0\right)}\right).
\end{equation}
Thanks to the $C^{0,1}$ regularity of the matrix $g$, there exists
some constant $C>0$ such that
\[
\left|F_{\sigma}^{g}-F_{\sigma}^{g\left(0\right)}\right|\leq C\sigma.
\] 
and, hence,  by \eqref{phi}, for  every $\sigma\geq 1$ ,
\begin{equation}\label{inter11}
|H(\sigma)-\epsilon \sigma^\epsilon|\leq C\sigma,
\end{equation}
where $C$ depends only on $\Lambda$ from \eqref{lipg}.
This also means that,  if $\sigma$ is small, $H(\sigma)$ is positive.

By \eqref{inter10} and \eqref{H}, we have
\begin{align}\label{eq:Formulation_Fwg_term}
\tau\gamma F_{w}^{g}\Delta_{g}\left(f^{2}\right)  =\tau\left(n-2\right)\gamma\phi\left(\sigma\right)\Delta_{g}\left(f^{2}\right)
  -2\tau\gamma H\left(\sigma\right)f\Delta_{g}f -2\tau\gamma H\left(\sigma\right)\left|\nabla_{g}f\right|^{2},
\end{align}
where we can write
\begin{equation}\label{inter12}
\gamma\phi\left(\sigma\right)\Delta_{g}\left(f^{2}\right)=\dive \left(\gamma\phi\left(\sigma\right)\nabla_g(f^2)\right)-\nabla_g(\gamma\phi\left(\sigma\right))\cdot\nabla_g(f^2) 
\end{equation}
almost everywhere in $B_1^+$. 

Let us estimate
\begin{eqnarray}\label{inter13}
\left|\nabla_g(\gamma\phi\left(\sigma\right))\cdot\nabla_g(f^2) \right|&\leq& %2
%\left|\nabla_g(\gamma\phi\left(\sigma\right)) \right| |f|\left|\nabla_g f \right|\\
%&\leq & 
2\left|\nabla_g\gamma \phi(\sigma)+\gamma\phi'\left(\sigma\right)\nabla_g\sigma \right||f|\left|\nabla_g f \right|\nonumber
\\&\leq& C\sigma^{\epsilon-1}|f|\left|\nabla_g f \right|
\end{eqnarray}
By \eqref{Pvtaus}, 
\begin{equation}\label{inter14}
\Delta_gf=P_{w,\tau}^{(s)}(f)-\tau^2 \frac{|\nabla_g w|^2}{w^2}f,
\end{equation}
and, hence, by \eqref{inter14} and \eqref{inter11}
\begin{align}\label{inter15}
-2\tau\gamma H\left(\sigma\right)f\Delta_{g}f&=-2\tau\gamma H\left(\sigma\right)fP_{w,\tau}^{(s)}(f)+2\tau^3\gamma\frac{|\nabla_g w|^2}{w^2} H\left(\sigma\right)f^2\geq\\
&\geq -2\tau\gamma (\epsilon \sigma^\epsilon+C\sigma)|f|\left|P_{w,\tau}^{(s)}(f)\right|+2\tau^3\gamma\frac{|\nabla_g w|^2}{w^2} (\epsilon \sigma^\epsilon-C\sigma)f^2\nonumber%\\
%&&\geq -2\tau\gamma (\epsilon \sigma^\epsilon+C\sigma)|f|\left|P_{w,\tau}^{(s)}(f)\right|+2C\tau^3\gamma\sigma^{-2} (\epsilon \sigma^\epsilon-C\sigma)f^2\nonumber\\
\end{align}
Now, by putting together \eqref{inter9},  \eqref{eq:Formulation_Fwg_term}, \eqref{inter12}, \eqref{inter13} and \eqref{inter14}, we finally get
\begin{align}\label{inter16}
\gamma&\frac{w^{2}}{\left|\nabla_{g}w\right|^{2}}\left(P_{w,\tau}f\right)^{2}\geq\gamma\frac{w^{2}}{\left|\nabla_{g}w\right|^{2}}\left(P_{w,\tau}^{\left(s\right)}f\right)^{2}+\nonumber\\
&+4\tau^{2}\gamma\frac{\left|\nabla_{g}w\right|^{2}}{w^{2}}\left(1-\frac{2(1-\delta)\epsilon\sigma^\epsilon}{\tau}\right)\left(A_{w}\left(f\right)\right)^{2}+\nonumber\\
&+\tau\left(4(1-\delta)\sigma^\epsilon\gamma\epsilon-C(1-\delta)\sigma -2\gamma H(\sigma)\vphantom{\frac{1}{1}}\right)\left|\nabla_{g}f\right|^{2}+\nonumber\\
&-C\tau
\sigma^{\epsilon-1}|f|\left|\nabla_{g}f\right|-2\tau\gamma(\epsilon\sigma^\epsilon+C\sigma)|f|\left|P_{w,\tau}^{(s)}(f)\right|+\nonumber\\
&+\tau^3\left(2\gamma \frac{|\nabla_gw|^2}{w^2}(\epsilon\sigma^\epsilon-C\sigma)-2\frac{\nabla_g w\cdot \nabla_g \gamma}{w}-8(1-\delta)C\tau^{-2}\gamma\epsilon\sigma^\epsilon\frac{|\nabla_gw|^2}{w^2}\right)f^2+\nonumber\\&+\dive\left(2\tau G_1 +\tau(n-2)\gamma\phi\left(\sigma\right)\nabla_g(f^2)\right)
\end{align}
By Young's inequality we can write
\begin{equation}\label{young1}
C\tau
\sigma^{\epsilon-1}|f|\left|\nabla_{g}f\right|\leq \frac{C}{\delta\gamma\epsilon}\tau\sigma^{\epsilon-2}f^2+\delta\tau\gamma\epsilon\sigma^\epsilon\left|\nabla_{g}f\right|^2
\end{equation}
and 
\begin{align}\label{young2}
2\tau\gamma(\epsilon\sigma^\epsilon+C\sigma)|f|\left|P_{w,\tau}^{(s)}(f)\right|\leq 
4\tau^2\gamma(\epsilon\sigma^\epsilon+C\sigma)^2\frac{|\nabla_g w|^2}{w^2}f^2+
\frac{1}{2}\frac{w^2}{|\nabla_g w|^2}\gamma \left|P_{w,\tau}^{(s)}(f)\right|^2
\end{align}
so that, finally, from \eqref{inter16}, \eqref{young1} and \eqref{young2},  we have
\begin{align}\label{inter17}
\gamma&\frac{w^{2}}{\left|\nabla_{g}w\right|^{2}}\left(P_{w,\tau}f\right)^{2}\geq\frac{\gamma}{2}\frac{w^{2}}{\left|\nabla_{g}w\right|^{2}}\left(P_{w,\tau}^{\left(s\right)}f\right)^{2}+\nonumber\\
&+4\tau^{2}\gamma\frac{\left|\nabla_{g}w\right|^{2}}{w^{2}}\left(1-\frac{2(1-\delta)\epsilon\sigma^\epsilon}{\tau}\right)\left(A_{w}\left(f\right)\right)^{2}+\nonumber\\
&+\tau\left(4(1-\delta)\sigma^\epsilon\gamma\epsilon-C(1-\delta)\sigma -2\gamma H(\sigma)-\delta\gamma\epsilon\sigma^\epsilon\vphantom{\frac{1}{1}}\right)\left|\nabla_{g}f\right|^{2}+\nonumber\\
&+\tau^3\left(2\gamma \frac{|\nabla_gw|^2}{w^2}(\epsilon\sigma^\epsilon-C\sigma)-2\frac{\nabla_g w\cdot \nabla_g \gamma}{w}-8(1-\delta)C\tau^{-2}\gamma\epsilon\sigma^\epsilon\frac{|\nabla_gw|^2}{w^2}\right.+\nonumber\\&\quad\left.-\frac{C}{\delta\gamma\epsilon}\tau^{-2}\sigma^{\epsilon-2}-4\tau^{-1}\gamma(\epsilon\sigma^\epsilon+C\sigma)^2\frac{|\nabla_gw|^2}{w^2}\right)f^2+\nonumber\\&+\dive\left(2\tau G_1 +\tau(n-2)\gamma\phi\left(\sigma\right)\nabla_g(f^2)\right)
\end{align}

Now we want to estimate from below all the coefficients appearing on the right hand side by selecting $\delta$ and choosing $\sigma$ small and $\tau$ big enough. 

Let us start from the coefficient of $|\nabla_gf^2|$. Let us choose $\delta=1/4$ and, by using \eqref{inter11} and recalling that $\sigma<1$, we can  write
\begin{align}\label{est1}
4(1-\delta)\sigma^\epsilon&\gamma\epsilon-C(1-\delta)\sigma -2\gamma H(\sigma)-\delta\gamma\epsilon\sigma^\epsilon
\geq\nonumber\\&\geq 4(1-\delta)\sigma^\epsilon\gamma\epsilon-C(1-\delta)\sigma -2\gamma \epsilon \sigma^\epsilon -2C\sigma\gamma -\delta\gamma\epsilon\sigma^\epsilon=\nonumber\\
&= 
 \frac{3}{4}\sigma^\epsilon\gamma\epsilon-\frac{3C}{4}\sigma  -2C\sigma\gamma \geq \frac{1}{2}\sigma^\epsilon\gamma\epsilon
\end{align}
for $\sigma \leq\overline{r}_1$ where $\overline{r}_1$ depends only on $\Lambda$,  $\lambda$,  $\epsilon$, $\gamma_0$ and $\|\gamma\|_{C^{0,1}}$. 

Now let us consider the coefficient of $f^2$ in \eqref{inter17} in which we take $\delta=1/4$.  Let us notice that, from \eqref{gradienteprec} and since $\epsilon\in(0,1)$, there exist $\overline{r}_2<1$ and  $\tau_0$   such that, for $\sigma<\overline{r}_2$ and $\tau>\tau_0$ we have
\begin{align}\label{est2}
2\gamma \frac{|\nabla_gw|^2}{w^2}&(\epsilon\sigma^\epsilon-C\sigma)-2\frac{\nabla_g w\cdot \nabla_g \gamma}{w}-8(1-\delta)C\tau^{-2}\gamma\epsilon\sigma^\epsilon\frac{|\nabla_gw|^2}{w^2}+\nonumber\\&-\frac{C}{\delta\gamma\epsilon}\tau^{-2}\sigma^{\epsilon-2}-4\tau^{-1}\gamma(\epsilon\sigma^\epsilon+C\sigma)^2\frac{|\nabla_gw|^2}{w^2}\geq\nonumber\\
& \geq \frac{\lambda^{-1}\gamma\sigma^{-2}}{2}(\epsilon\sigma^\epsilon-C\sigma)-C\sigma^{-1}-6C\tau^{-2}\gamma\epsilon\sigma^{\epsilon-2}+\nonumber\\&-\frac{4C}{\gamma\epsilon}\tau^{-2}\sigma^{\epsilon-2}-4C\tau^{-1}\gamma(\epsilon\sigma^\epsilon+C\sigma)^2\sigma^{-2}\geq\frac{ \gamma \lambda^{-1}\epsilon\sigma^{\epsilon-2}}{4}
\end{align}
Hence, finally, by \eqref{inter17}, \eqref{est1} and \eqref{est2}, we get
\begin{align}\label{inter18}
\gamma\frac{w^{2}}{\left|\nabla_{g}w\right|^{2}}\left(P_{w,\tau}f\right)^{2}&\geq\frac{\gamma}{2}\frac{w^{2}}{\left|\nabla_{g}w\right|^{2}}\left(P_{w,\tau}^{\left(s\right)}f\right)^{2}+2\tau^{2}\gamma\frac{\left|\nabla_{g}w\right|^{2}}{w^{2}}\left(A_{w}\left(f\right)\right)^{2}+\nonumber\\
&+\frac{1}{2}\tau\sigma^\epsilon\gamma\epsilon\left|\nabla_{g}f\right|^{2}+\tau^3f^2\frac{ \gamma \lambda^{-1}\epsilon\sigma^{\epsilon-2}}{4}+\nonumber\\
&+\dive\left(2\tau G_1 +\tau(n-2)\gamma\phi\left(\sigma\right)\nabla_g(f^2)\right)
\end{align}
 for $\tau\geq \tau_0$ and $\sigma\leq \overline{r}:=\min\{\overline{r}_0, \overline{r}_1, \overline{r}_2\}$, that is for $x\in B^+_{\overline{r}}$.
 By simple calculations \eqref{inter18} gives \eqref{senzaint} and \eqref{G2}.
 \end{proof}

Proposition \ref{propsenzaint} holds also in $B_1^-$, where $f$ denotes $f^-$ and $\gamma$ is any positive function in $C^{0,1}\left(B_1^-\right)$.

\subsection{Integral estimates for $f$}
Now we integrate the obtained pointwise estimate \eqref{senzaint} on $B_1^+$ and the corresponding estimate on $B_1^-$ and integrate by parts the last term $\dive G_2$. Since estimate \eqref{senzaint} holds in the smaller set  $B_{\overline{r}}$, we will restrict to functions  $f\in \mathcal{H}(B_1)$ with support  contained in $B_{\overline{r}}\setminus\{0\}$.  For this reason, this integration by parts does not give any contribution on $\partial B_1$.

In order to obtain a typical Carleman estimate, we  neglect the first positive terms on the right hand side of \eqref{senzaint}.

\begin{proposition}\label{propconint}
Let $g$ satisfy assumptions \eqref{ellipticity}, \eqref{lipg},  \eqref{ginzero}, and \eqref{gatzero}, let $\gamma^+\in C^{0,1}(B_1^+)$ and $\gamma^-\in C^{0,1}(B_1^-)$ such that $\gamma^\pm\geq \gamma_0>0$ in $B_1^\pm$ and let $w$ as in \eqref{w}.
Then, there exist $\overline{r}\in(0,1)$ and $\tau_0>0$, depending only on $\lambda$, $\Lambda$,  $\gamma_0$ and \,  $\|\gamma^\pm\|_{C^{0,1}(B_1^\pm)}$ and $\epsilon$ such that, for any $f\in \mathcal{H}(B_1)$ supported in $B_{\overline{r}}\setminus\{0\}$
and for any $\tau\geq \tau_0$ we have
\begin{align}\label{conint}
\intop_{B_1}&\frac{\gamma w^{2}}{\left|\nabla_{g}w\right|^{2}}\left(P_{w,\tau}f\right)^{2}\geq2\tau^{2}\intop_{B_1}\frac{\gamma\left|\nabla_{g}w\right|^{2}}{w^{2}}\left(A_{w}\left(f\right)\right)^{2}+\\
&+\frac{\tau\epsilon}{2}\intop_{B_1}\gamma\sigma^\epsilon\left|\nabla_{g}f\right|^{2}+\frac{\tau^3\epsilon\lambda^{-1}}{4}\intop_{B_1} \gamma \sigma^{\epsilon-2}f^2+\nonumber\\&+\intop_{\{x_n=0\}}\left[\left(4\tau w\left( \frac{\nabla_{g}w\cdot\nabla_{g}f}{|\nabla_{g}w|^2}\right) +2\tau(n-2)\phi\left(\sigma\right)f\vphantom{\frac{l}{l}}\right)\left(\gamma\nabla_{g} f\cdot \nu_g\right)\right],\nonumber
\end{align}
where $\gamma=\gamma^+\chi_{B_1^+}+\gamma^-\chi_{B_1^-}$. %Notice that $\gamma\in C^{0,1}\left(B_1\setminus \{x_n=0\}\right)$, but is is not continuous on $\{x_n=0\}$.
\end{proposition}
\begin{proof}
We start by integrating on $B_1^+$ inequality \eqref{senzaint} (where we neglect the positive term containing $P_{w,\tau}^{(s)}f$ and take $\gamma=\gamma^+$) and get
\begin{align}\label{f1}
\intop_{B_1^+}\gamma^+\frac{w^{2}}{\left|\nabla_{g}w\right|^{2}}\left(P_{w,\tau}f\right)^{2}\geq&2\tau^{2}\intop_{B_1^+}\gamma^+\frac{\left|\nabla_{g}w\right|^{2}}{w^{2}}\left(A_{w}\left(f\right)\right)^{2}+\frac{\tau}{2}\intop_{B_1^x}\sigma^\epsilon\gamma^+\epsilon\left|\nabla_{g}f\right|^{2}+\nonumber\\
&+\frac{\tau^3}{4}\intop_{B_1^+} \gamma^+ \lambda^{-1}\epsilon\sigma^{\epsilon-2}f^2+\intop_{B_1^+}\dive G_2
\end{align}
Now let us integrate by parts the last term in this inequality.  Let us recall that $\nu=-e_n$ is the outgoing normal with respect to $B_1^+$ on the interface $\{x_n=0\}$.
Let us also recall that the support of $f$ is far from $\partial B_1$, so that
\begin{equation*}%\label{f2}
\intop_{B_1^+}\dive G_2=\intop_{\{x_n=0\}}G_2\cdot \nu_g
\end{equation*}
Now, by using \eqref{Bv}, we calculate
\begin{align*}%\label{f3}
G_2\cdot \nu_g &=\left(2\tau^3 \frac{ f^2}{w}-2\tau \left|\nabla_{g}f\right|^{2}\frac{w}{|\nabla_gw|^2}\right)\gamma^+\nabla_gw\cdot \nu_g+\nonumber\\&+\left(4\tau\left( B_w\cdot\nabla_{g}f\right) +2\tau(n-2)\phi\left(\sigma\right)f\vphantom{\frac{l}{l}}\right)\gamma^+\nabla_{g}f \cdot\nu_g
\end{align*}
By \eqref{w}  and \eqref{sigma} and since $\nu=-e_n$  we have, for $x=(x^\prime,0)$,
\begin{equation*}
\nabla_gw\cdot\nu_g=\psi^\prime(\sigma)\nabla_g\sigma\cdot\nu_g= -\psi^\prime(\sigma)\sum_{j=1}^{n-1}g^{nj}(x^\prime,0)\frac{x_j}{|x|}.
\end{equation*}
By assumption \eqref{gatzero}, we get
\begin{equation}\label{f4}
\nabla_gw\cdot \nu_g=0.
\end{equation}
Then, we have, on $\{x_n=0\}$,
\begin{equation}\label{f4bis}
G_2\cdot\nu_g=\left(4\tau w\left( \frac{\nabla_{g}w\cdot\nabla_{g}f}{|\nabla_{g}w|^2}\right) +2\tau(n-2)\phi\left(\sigma\right)f\vphantom{\frac{l}{l}}\right)\gamma^+\nabla_{g}f \cdot \nu_g.
\end{equation}
So, finally, from \eqref{f1} and \eqref{f4bis},
\begin{align}\label{f5piu}
\intop_{B_1^+}\gamma^+&\frac{w^{2}}{\left|\nabla_{g}w\right|^{2}}\left(P_{w,\tau}f\right)^{2}\geq2\tau^{2}\intop_{B_1^+}\gamma^+\frac{\left|\nabla_{g}w\right|^{2}}{w^{2}}\left(A_{w}\left(f\right)\right)^{2}+\nonumber\\&+\frac{\tau}{2}\intop_{B_1^+}\sigma^\epsilon\gamma^+\epsilon\left|\nabla_{g}f\right|^{2}+\frac{\tau^3}{4}\intop_{B_1^+} \gamma^+ \lambda^{-1}\epsilon\sigma^{\epsilon-2}f^2+\\&+\intop_{\{x_n=0\}}\!\!\left(4\tau w\left( \frac{\nabla_{g}w\cdot\nabla_{g}f^+}{|\nabla_{g}w|^2}\right) +2\tau(n-2)\phi\left(\sigma\right)f^+\right)\gamma^+\nabla_{g}f^+ \cdot \nu_g,\nonumber
\end{align}
Now, let us do the same in $B_1^-$, by using $\gamma=\gamma^-$. The unit normal at $\{x_n=0\}$ which is an exterior normal with respect to $B_1^-$ it is now $-\nu=e_n$ so that
\begin{align}\label{f5meno}
\intop_{B_1^-}\gamma^-\frac{w^{2}}{\left|\nabla_{g}w\right|^{2}}&\left(P_{w,\tau}f\right)^{2}\geq2\tau^{2}\intop_{B_1^-}\gamma^-\frac{\left|\nabla_{g}w\right|^{2}}{w^{2}}\left(A_{w}\left(f\right)\right)^{2}+\nonumber\\&+\frac{\tau}{2}\intop_{B_1^-}\sigma^\epsilon\gamma^-\epsilon\left|\nabla_{g}f\right|^{2}+\frac{\tau^3}{4}\intop_{B_1^-} \gamma^- \lambda^{-1}\epsilon\sigma^{\epsilon-2}f^2+\\&-\intop_{\{x_n=0\}}\!\!\left(4\tau w\left( \frac{\nabla_{g}w\cdot\nabla_{g}f^-}{|\nabla_{g}w|^2}\right) +2\tau(n-2)\phi\left(\sigma\right)f^-\right)\gamma^-\nabla_{g}f^- \cdot \nu_g.\nonumber
\end{align}
By adding \eqref{f5piu} and \eqref{f5meno} 
 we finally have \eqref{conint}.
\end{proof}

\subsection{Proof of Theorem \ref{teopiatta}}
\begin{proof}
Let us set 
\begin{equation*}
f=w^{-\tau}u.
\end{equation*}
We have
\begin{equation}\label{Pfeu}
P_{w,\tau}\left(f\right)=w^{-\tau}\Delta_{g}u,
\end{equation}
and
\begin{equation}\label{gradfeu}
\nabla_g f=-\tau w^{-\tau-1}u\nabla_g w+w^{-\tau}\nabla_g u.
\end{equation}
Notice that, if $u\in \mathcal{H}(B_1)$ and is supported in $B_{\overline{r}}^+\setminus\{0\}$, the same holds for $f$.

Let us consider the left hand side of \eqref{conint}. For sake of shortness we only consider the integrals on the upper half ball.
We have, by \eqref{Pfeu} and \eqref{gradienteprec},
\begin{equation}\label{LHS}
\intop_{B_1^+}\frac{\gamma w^{2}}{\left|\nabla_{g}w\right|^{2}}\left(P_{w,\tau}f\right)^{2}=\intop_{B_1^+}\frac{\gamma w^{2}}{\left|\nabla_{g}w\right|^{2}}w^{-2\tau}\left(\Delta_{g}u\right)^2\leq 4\lambda\intop_{B_1^+}
\gamma\sigma^2w^{-2\tau}\left(\Delta_{g}u\right)^2
\end{equation}
On the right hand side, we will omit the positive termi containing $\left(A_{w}\left(f\right)\right)^{2}$. 
For the second and third term of the right hand side,  from \eqref{gradfeu} we write, for some positive $\delta<1$ (that we will choose later)
 \begin{equation*}
\left|\nabla_g f\right|^2\geq (1-\delta)w^{-2\tau}\left|\nabla_g u\right|^2-(\delta^{-1}-1)\tau^2 w^{-2\tau-2}u^2\left|\nabla_g w\right|^2,
\end{equation*}
which, by \eqref{gradienteprec}, gives
 \begin{equation*}%\label{gradfeu1}
\left|\nabla_g f\right|^2\geq (1-\delta)w^{-2\tau}\left|\nabla_g u\right|^2-(\delta^{-1}-1)\lambda\tau^2 w^{-2\tau}\sigma^{-2}u^2,
\end{equation*}
and, hence,
 \begin{align*}%\label{gradfeu2}
&\frac{\tau\epsilon}{2}\intop_{B_1^+}\gamma\sigma^\epsilon\left|\nabla_{g}f\right|^{2}+\frac{\tau^3\epsilon\lambda^{-1}}{4}\intop_{B_1^+} \gamma \sigma^{\epsilon-2}f^2\geq\\
&\geq  (1-\delta)
\frac{\tau\epsilon}{2}\intop_{B_1^+}\gamma\sigma^\epsilon w^{-2\tau}\left|\nabla_{g}u\right|^{2}+\tau^3\epsilon\left(\frac{\lambda^{-1}}{4}-\frac{(\delta^{-1}-1)\lambda}{2}\right)\intop_{B_1^+} \gamma \sigma^{\epsilon-2}w^{-2\tau}u^2\nonumber
\end{align*}
By choosing $\delta^{-1}=1+\frac{\lambda^{-2}}{8}$, we get
 \begin{align*}%\label{gradfeu3}
&\frac{\tau\epsilon}{2}\intop_{B_1^+}\gamma\sigma^\epsilon\left|\nabla_{g}f\right|^{2}+\frac{\tau^3\epsilon\lambda^{-1}}{4}\intop_{B_1^+} \gamma \sigma^{\epsilon-2}f^2\geq\\
&\geq \frac{}{}
\frac{\tau\epsilon\lambda^{-2}}{2(8+\lambda^{-2})}\intop_{B_1^+}\gamma\sigma^\epsilon w^{-2\tau}\left|\nabla_{g}u\right|^{2}+\frac{\tau^3\epsilon\lambda^{-1}}{8}
\intop_{B_1^+} \gamma \sigma^{\epsilon-2}w^{-2\tau}u^2\nonumber
\end{align*}
For the boundary terms, we recall that, by assumption \eqref{gatzero},  \eqref{f4} holds, hence, on $\{x_n=0\}$
\begin{equation*}
\gamma\nabla_gf\cdot\nu_g=\gamma w^{-\tau}\nabla_g u\cdot \nu_g
\end{equation*}
and
\begin{equation*}
\nabla_gw\cdot\nabla_gf= -\tau w^{-\tau-1}u\left|\nabla_g w\right|^2+w^{-\tau}\nabla_gw\cdot\nabla_gu
\end{equation*}
hence
\begin{align}\label{bdu}
&\intop_{\{x_n=0\}}\!\!\left(4\tau w\left( \frac{\nabla_{g}w\cdot\nabla_{g}f^+}{|\nabla_{g}w|^2}\right) +2\tau(n-2)\phi\left(\sigma\right)f^+\right)\gamma^+\nabla_{g}f^+ \cdot \nu_g=\\
&\intop_{\{x_n=0\}}\!\!w^{-2\tau}\left(4\tau w\left( \frac{\nabla_{g}w\cdot\nabla_{g}u^+}{|\nabla_{g}w|^2}\right) +\left(2\tau(n-2)\phi\left(\sigma\right)-4\tau^2\right)u^+\vphantom{\frac{l}{l}}\right)\gamma^+\nabla_{g}u^+ \cdot \nu_g.\nonumber
\end{align}
From \eqref{conint},  \eqref{LHS},  \eqref{gradfeu} and \eqref{bdu}, we finally get \eqref{uconint}.
\end{proof}

\section{Doubling inequality}
We now want to prove a doubling inequality for solutions of some elliptic partial differential equation with piecewise regular coefficients. Let us start with a slight modification of the Carleman estimate we proved int the previous section.

\begin{lemma}\label{lemmaconintdoub}
Let $g$ satisfy assumptions \eqref{ellipticity}, \eqref{lipg},  \eqref{ginzero}, and \eqref{gatzero}, let $\gamma^+\in C^{0,1}(B_1^+)$ and $\gamma^-\in C^{0,1}(B_1^-)$ such that $\gamma^\pm\geq \gamma_0>0$ in $B_1^\pm$ and let $w$ as in \eqref{w}.
Then, there exist $\overline{r}_1<\overline{r}$ of Proposition \eqref{propconint} and $C$, depending only on $\lambda$, $\Lambda$,  $\gamma_0$,  $\|\gamma^\pm\|_{C^{0,1}(B_1^\pm)}$ and $\epsilon$ such that, for any $0<r<\overline{r}_1$ and for every $f\in \mathcal{H}(B_1)$ supported in $B_{\overline{r}_1}\setminus B_r$
 we have
\begin{equation}\label{fsdoub}
r\intop_{B_1}
\gamma\frac{\left|\nabla_{g}w\right|^{2}}{w^{3}}f^2\leq C
\intop_{B_1}\frac{\gamma\left|\nabla_{g}w\right|^{2}}{w^{2}}\left(A_{w}\left(f\right)\right)^{2}
\end{equation}
where $\gamma=\gamma^+\chi_{B_1^+}+\gamma^-\chi_{B_1^-}$. 
\end{lemma}
\begin{proof}
By \eqref{Av} and \eqref{Fgv}, we can write
\begin{equation}\label{fsdoub1}
\frac{\gamma\left|\nabla_{g}w\right|^{2}}{w^{3}}A_{w}\left(f\right)f=
\frac{1}{2}\dive\left(\frac{\gamma f^2}{w^2}\nabla_gw\right)+\frac{1}{2}\frac{\gamma\left|\nabla_{g}w\right|^{2}}{w^3}f^2\left(1-\frac{w}{\gamma}\frac{\nabla_g\gamma\cdot\nabla_gw}{\left|\nabla_{g}w\right|^{2}}\right).
\end{equation}
Now,  if we integrate \eqref{fsdoub1} on  $B^+_{\overline{r}_1}\setminus B^+_r$  and we recall that $f=0$ on $\partial B_{\overline{r}_1}\cup\partial B_r$ and $\nabla_g w\cdot \nu=0$ on $\{x_n=0\}$, we get
\begin{equation*}%\label{fsdoub2}
\intop_{B^+_{\overline{r}_1}\setminus B^+_r}\!\!\!\!\frac{\gamma\left|\nabla_{g}w\right|^{2}}{w^{3}}A_{w}\left(f\right)f=\frac{1}{2}\intop_{B^+_{\overline{r}_1}\setminus B^+_r}\!\!\!\!\frac{\gamma\left|\nabla_{g}w\right|^{2}}{w^3}f^2\left(1-\frac{w}{\gamma}\frac{\nabla_g\gamma\cdot\nabla_gw}{\left|\nabla_{g}w\right|^{2}}\right).
\end{equation*}
The same holds integrating on $B^-_{\overline{r}_1}\setminus B^-_r$.

By the assumption on $\gamma$ and by  \eqref{gradienteprec}, we have
\begin{equation*}
1-\frac{w}{\gamma}\frac{\nabla_g\gamma\cdot\nabla_gw}{\left|\nabla_{g}w\right|^{2}}\geq 1-C\sigma
\end{equation*}
where $C$ depends on $\lambda$, $\gamma_0$ and $\|\gamma\|_{C^{0,1}}$.
Hence, if we take $\overline{r}_1<\overline{r}$ such that $C\overline{r}_1\leq \frac{1}{2}$, we have
\begin{equation}\label{fsdoub3}
\intop_{B^+_{\overline{r}_1}\setminus B^+_r}\!\!\!\!\frac{\gamma\left|\nabla_{g}w\right|^{2}}{w^{3}}A_{w}\left(f\right)f\geq \frac{1}{4}
\intop_{B^+_{\overline{r}_1}\setminus B^+_r}\!\!\!\!\frac{\gamma\left|\nabla_{g}w\right|^{2}}{w^3}f^2.
\end{equation}
Now, by Cauchy inequality, and since we take $\sigma\geq r$ (and, hence $w\geq \Psi(r)$),  we can write,
\begin{align}\label{fsdoub4}
\intop_{B^+_{\overline{r}_1}\setminus B^+_r}&\!\!\!\!\frac{\gamma\left|\nabla_{g}w\right|^{2}}{w^{3}}A_{w}\left(f\right)f\leq
\left(\intop_{B^+_{\overline{r}_1}\setminus B^+_r}\!\!\!\!\frac{\gamma\left|\nabla_{g}w\right|^{2}}{w^{2}}\left(A_{w}\left(f\right)\right)^2\right)^{\frac{1}{2}}\!\!
\left(\intop_{B^+_{\overline{r}_1}\setminus B^+_r}\!\!\!\!\frac{\gamma\left|\nabla_{g}w\right|^{2}}{w^{4}}f^2\right)^{\frac{1}{2}}\!\!\leq\nonumber
\\&\leq \frac{4}{\sqrt{\Psi(r)}}\left(\intop_{B^+_{\overline{r}_1}\setminus B^+_r}\!\!\!\!\frac{\gamma\left|\nabla_{g}w\right|^{2}}{w^{2}}\left(A_{w}\left(f\right)\right)^2\right)^{\frac{1}{2}}
\left(\intop_{B^+_{\overline{r}_1}\setminus B^+_r}\!\!\!\!\frac{\gamma\left|\nabla_{g}w\right|^{2}}{w^{3}}f^2\right)^{\frac{1}{2}}.
\end{align}
From \eqref{fsdoub3} and \eqref{fsdoub4} (and repeating the same procedure on $B^-$ we finally have 
\begin{equation*}%\label{fsdoub5}
\psi(r)\intop_{B_{\overline{r}_1}\setminus B_r}\frac{\gamma\left|\nabla_{g}w\right|^{2}}{w^{3}}f^2\leq 16 \intop_{B_{\overline{r}_1}\setminus B_r}\frac{\gamma\left|\nabla_{g}w\right|^{2}}{w^{2}}\left(A_{w}\left(f\right)\right)^2.
\end{equation*}
By \eqref{asintpsi},  for $\overline{r}_1$ small enough, we have \eqref{fsdoub} for every $r$ such that  $0<r<\overline{r}_1$.
\end{proof}
Now, let us insert \eqref{fsdoub} into Proposition \ref{propsenzaint} and write the analogue of Theorem \ref{teopiatta}.
\begin{proposition}\label{propdoub}
Let $g$ satisfy assumptions \eqref{ellipticity}, \eqref{lipg},  \eqref{ginzero}, and \eqref{gatzero}, let $\gamma^+\in C^{0,1}(B_1^+)$ and $\gamma^-\in C^{0,1}(B_1^-)$ such that $\gamma^\pm\geq \gamma_0>0$ in $B_1^\pm$ and let $w$ as in \eqref{w}.
Then, there exist $\overline{r}_1\in(0,1)$, $C>0$ and $\tau_0>0$, depending only on $\lambda$, $\Lambda$,  $\gamma_0$,  $\|\gamma^\pm\|_{C^{0,1}(B_1^\pm)}$ and $\epsilon$ such that,  for any $0<r<\overline{r}_1$ and for every $u\in \mathcal{H}(B_1)$ supported in $B_{\overline{r}_1}\setminus B_r$
and for any $\tau\geq \tau_0$ we have
\begin{align}\label{uconintdoub}
&4\lambda\intop_{B_1}\gamma\sigma^2 w^{-2\tau}\left(\Delta_gu\right)^2 
\geq \frac{\tau\epsilon}{2(1+8\lambda^2)}\intop_{B_1}\gamma\sigma^\epsilon w^{-2\tau}\left|\nabla_{g}u\right|^{2}+\\
&+\frac{\tau^3\epsilon}{8\lambda}\intop_{B_1} \gamma \sigma^{\epsilon-2}w^{-2\tau}u^2+\frac{\tau^2r}{C}\intop_{B_1}\sigma^{-2}w^{-2\tau-1}u^2+
\nonumber\\
&+\intop_{\{x_n=0\}}\!\!w^{-2\tau}\left[\left(4\tau w\left( \frac{\nabla_{g}w\cdot\nabla_{g}u}{|\nabla_{g}w|^2}\right) +\left(2\tau(n-2)\phi\left(\sigma\right)-4\tau^2\right)u\vphantom{\frac{l}{l}}\right)\left(\gamma\nabla_{g} u\cdot \nu_g\right)\right],\nonumber
\end{align}
where $\gamma=\gamma^+\chi_{B_1^+}+\gamma^-\chi_{B_1^-}$. 
\end{proposition}
\begin{remark}\label{remarkreg}
Notice that, by a density argument, inequality $\eqref{uconintdoub}$ holds for every $u\in H^1(B_1)$ supported in $B_{\overline{r}_1}\setminus B_r$ such that $u^+\in H^2(B_1^+)$ and $u^-\in H^2(B_1^-)$.
\end{remark}
\begin{proof}
    Let us take $\overline{r}_1$ as in Lemma \ref{lemmaconintdoub}, $0<r<\overline{r}_1$ and $f=w^{-\tau}u$ that will be supported in $B_{\overline{r}_1}\setminus B_r$.
   By putting together 
    \eqref{conint} and  \eqref{fsdoub}, we get that
\begin{align*}%\label{doubforf}
\intop_{B_1}&\frac{\gamma w^{2}}{\left|\nabla_{g}w\right|^{2}}\left(P_{w,\tau}f\right)^{2}\geq\frac{2\tau^{2}r}{C}\intop_{B_1}\gamma\frac{\left|\nabla_{g}w\right|^{2}}{w^{3}}f^2+\\
&+\frac{\tau\epsilon}{2}\intop_{B_1}\gamma\sigma^\epsilon\left|\nabla_{g}f\right|^{2}+\frac{\tau^3\epsilon\lambda^{-1}}{4}\intop_{B_1} \gamma \sigma^{\epsilon-2}f^2+\nonumber\\&+\intop_{\{x_n=0\}}\left[\left(4\tau w\left( \frac{\nabla_{g}w\cdot\nabla_{g}f}{|\nabla_{g}w|^2}\right) +2\tau(n-2)\phi\left(\sigma\right)f\vphantom{\frac{l}{l}}\right)\left(\gamma\nabla_{g} f\cdot \nu_g\right)\right],\nonumber
\end{align*}
    where $C$ is given in Lemma \ref{lemmaconintdoub}.

Now, following the same procedure of the proof of Theorem 2.1, and using \eqref{gradienteprec} we get \eqref{uconintdoub}.
\end{proof}

Let $\tilde{a}^{\pm}\in C^{0,1}\left(\mathbb{R}^{d}\right)$ be positive functions, and let us set
\begin{equation*}
    \tilde{a}=\begin{cases}
\tilde{a}^{+} & \text{in }B_{1}^{+}\\
\tilde{a}^{-} & \text{in }B_{1}^{-}
\end{cases}
\end{equation*}
Let also assume that 
\begin{equation}\label{assa}
    \tilde{a}^\pm\geq \gamma_0>0 \mbox{ in }B_1^\pm. 
\end{equation}

Consider the following second order elliptic equation with piecewise regular 
coefficients: 
\begin{equation}\label{eq:piecewiseeq_pb}
\dive \left(\tilde{a}\nabla_{g}U\right)=0.
\end{equation}
A solution of this equation is a function $U\in H^1(\Omega)$ such that, if we set
\begin{equation*}
    U=\begin{cases}
U^{+} & \text{in }B_{1}^{+}\\
U^{-} & \text{in }B_{1}^{-}
\end{cases}
\end{equation*}
we have that 
\begin{equation*}
    \dive\left(\tilde{a}^\pm\nabla_g U^\pm\right)=0 \mbox{ in the interior of } B_1^\pm
\end{equation*}
and the transmission conditions
\begin{equation}\label{transm}
\begin{cases}
\tilde{a}^{+}\nabla_{g}U^{+}\cdot e_{n}=\tilde{a}^{-}\nabla_{g}U^{-}\cdot e_{n} & \text{on }B_{1}\cap\left\{ x_{n}=0\right\}, \\
U^{+}=U^{-} & \text{on }B_{1}\cap\left\{ x_{n}=0\right\} 
\end{cases}
\end{equation}
are satisfied.

Our goal is to build  a doubling inequality
for solutions to  \eqref{eq:piecewiseeq_pb}. 

Let us first show this preliminary result:
\begin{lemma}
\label{lem:one_lem_before_doubling}Let $U$ be a solution of \eqref{eq:piecewiseeq_pb}. 
There exist $C>0$,  $\tau_2$ and $\overline{r}_1$ depending only on $\lambda$, $\Lambda$, $\gamma_0$ and $\|\tilde{a}^\pm\|_{C^{0,1}(B_1^\pm)}$, such that
\begin{equation}\label{tesilemma}
    \frac{1}{R^{1+2\tau}} \intop_{B_{R}} U^{2} + \frac{R}{\left(4r\right)^{2+2\tau}}\intop_{B_{4r}}U^{2}
\leq C\left(\frac{1}{r^{2+2\tau}}\intop_{B_{2r}}U^{2}+\frac{1}{\overline{r}_{1}^{2+2\tau}}\intop_{B_{\overline{r}_{1}}}U^{2}\right)
\end{equation}
%\begin{equation}%\label{tesilemma}
 %   R^{-1-2\tau} \intop_{B_{R}} U^{2} + R\left(4r\right)^{-2-2\tau}\intop_{B_{4r}}U^{2}
%\leq C\left(r^{-2-2\tau}\intop_{B_{2r}}U^{2}+\overline{r}_{1}^{-2-2\tau}\intop_{B_{\overline{r}_{1}}}U^{2}\right)
%\end{equation}
for every  $\tau>\tau_{2}$ and every $0<4r<R\leq\frac{\overline{r}_{1}}{4}$.
\end{lemma}
\begin{proof}
We want to apply Proposition \ref{propdoub} to a solution $U$ to \eqref{eq:piecewiseeq_pb}. For this, we choose $\gamma^\pm=\tilde{a}^\pm$.
In order to fulfill the assumption on the support of the function, we need to multiply $U$
with some cutoff function so that its support set is contained in
the ring $B_{\overline{r}_{1}}\setminus B_{r}$.

Let $\overline{r}_1$ be as in Proposition \ref{propdoub} and let us fix $r$ and $R$ such that
$0<4r<R\leq\frac{\overline{r}_{1}}{4}$.

Let us take a cutoff function $\eta\in C^{\infty}\left(\left[0,+\infty\right)\right)$
satisfying
\begin{enumerate}
\item $0\leq\eta\leq1$
\item $\eta=\begin{cases}
0 & \text{in }\left[0,r\right)\cup\left(\frac{\overline{r}_{1}}{2},1\right)\\
1 & \text{in }\left(2r,\frac{\overline{r}_{1}}{4}\right)
\end{cases}$
%$\eta=\begin{cases}
%0 & \text{in }\left[0,r\right)\cup\left(\frac{2}%{3}\overline{r}_{1},+\infty\right)\\
%1 & \text{in }\left(2r,\frac{\overline{r}_{1}}{2}\right)
%\end{cases}$
\item For $0\leq k\leq2$, $\left|\frac{d^{k}\eta}{dt^{k}}\left(t\right)\right|\leq\begin{cases}
Cr^{-k} & \text{in}\left(r,2r\right)\\
C\overline{r}_{1}^{-k} & \text{in}\left(\frac{\overline{r}_{1}}{4},\frac{\overline{r}_{1}}{2}\right)%\text{in}\left(\frac{\overline{r}_{1}}{2},\frac{2}{3}\overline{r}_{1}\right)
\end{cases}$
\end{enumerate}
Let us define $\xi\left(x\right)=\eta\left(\left|x\right|\right)$. 

Since $U\in H^{1}\left(B_{1}\right)$ solves \eqref{eq:piecewiseeq_pb}, by the transmission conditions \eqref{transm} we have that $U$ is continuous on $\{x_n=0\}$. Also, by regularity results for solutions to \eqref{eq:piecewiseeq_pb}, we have that $U^\pm\in H^2(B^\pm_1)$ (\cite{La-Ur}). 

This implies that  $\xi U^\pm\in H^2\left(B_{1}^\pm\right)$ and, by the definition of $\xi$ it is  supported in $B_{\overline{r}_{1}}\setminus B_{r}$. Thus, by Remark \ref{remarkreg} we can apply  Proposition \ref{propdoub}
to $\xi U$.

Let us start from boundary terms. Since $\xi U$ is continuous on the interface $\{x_n=0\}$, also its tangential derivatives are continuous on the interface. 
For this reason, since by assumption \eqref{gatzero},
\begin{equation*}
\nabla_{g}w\cdot\nabla_{g}(\xi U)=\sum_{i,j=1}^{n-1}g^{ij}\partial_{j}w\partial_{i}(\xi U)
\end{equation*}
this is also continuous through $\{x_n=0\}$.

This implies that, on the interface $\{x_n=0\}$ we have
\begin{align}\label{bdxiu1}
&\left[\left(4\tau w\left( \frac{\nabla_{g}w\cdot\nabla_{g}(\xi U)}{|\nabla_{g}w|^2}\right) +\left(2\tau(n-2)\phi\left(\sigma\right)-4\tau^2\right)(\xi U)\vphantom{\frac{l}{l}}\right)\left( \tilde{a}\nabla_{g} (\xi U)\cdot \nu_g\right)\right]=\nonumber\\
&\left(4\tau w\left( \frac{\nabla_{g}w\cdot\nabla_{g}(\xi U)}{|\nabla_{g}w|^2}\right) +\left(2\tau(n-2)\phi\left(\sigma\right)-4\tau^2\right)(\xi U)\right)\left[\tilde{a} \nabla_{g} (\xi U)\cdot \nu_g\right]
\end{align}

Notice now that the boundary term in \eqref{bdxiu1} is zero because
\begin{equation}\label{bdxiu}
    \left[\tilde{a} \nabla_{g} (\xi U)\cdot \nu_g\right]=
\xi\left[\tilde{a} \nabla_{g} U\cdot \nu_g\right]+U\left[\tilde{a}\nabla_{g} \xi\cdot \nu_g\right]=0
\end{equation}
and the first term is zero for the transmission conditions \eqref{transm} %(remember that $\gamma=a$) 
while the second is zero by \eqref{gatzero} and because  $\xi$ is a radial function (the same reasons why \eqref{f4} holds).

We now use assumption \eqref{assa} and boundedness of $\tilde{a}^\pm$ and estimates \eqref{stimewl} to get, 
From \eqref{uconintdoub} (with $\epsilon=1$) and \eqref{bdxiu}, that there exist a constant $C$ (depending only on $\lambda$, $\Lambda$, $\gamma_0$ and $\|\tilde{a}^\pm\|_{C^{0,1}(B_1^\pm)}$) such  that
\begin{align}\label{4.19}
&C\intop_{B_1} w^{2-2\tau}\left(\Delta_g(\xi U)\right)^2 
\geq \tau\intop_{B_1} w^{1-2\tau}\left|\nabla_{g}(\xi U)\right|^{2}+\nonumber\\
&+\tau^3\intop_{B_1}w^{-1-2\tau}(\xi U)^2+\tau^2r\intop_{B_1}w^{-2\tau-3}(\xi U)^2
\end{align}
Let us now estimate from above  the left hand side in \eqref{4.19}. We first have
\begin{equation*}
 \left(\Delta_{g}\left(\xi U\right)\right)^{2}\leq2\xi^{2}\left(\Delta_{g}U\right)^{2}+16\left(U^{2}\left|\Delta_{g}\xi\right|^{2}+\left|\nabla_{g}U\right|^{2}\left|\nabla_{g}\xi\right|^{2}\right), 
\end{equation*}
and by using \eqref{stimewl} and the equation in \eqref{eq:piecewiseeq_pb}, that gives
\begin{equation*}
    \Delta_{g} U=-\frac{\nabla_{g} \tilde{a}}{\tilde{a}}\nabla_{g} U
    \end{equation*}
we have
%\begin{align}
 %   \sum_\pm\intop_{B_1^\pm}   w^{2-2\tau}\left(\Delta_g(\xi U)\right)^2 &\leq
  %  C  \sum_\pm\intop_{B_{\frac{\overline{r}_{1}}%{4}}^{\pm}\setminus B_{2r}^{\pm}}w^{2-2\tau}\left|\nabla_{g}U\right|^{2}  \\&+C\left(T_{1}+T_{2}\right)
%\end{align}
\begin{equation}\label{4.21}
  \intop_{B_1}   w^{2-2\tau}\left(\Delta_g(\xi U)\right)^2 \leq
    C \left( \intop_{B_{\frac{\overline{r}_{1}}{4}}\setminus B_{2r}}w^{2-2\tau}\left|\nabla_{g}U\right|^{2}+T_{1}+T_{2}\right)
\end{equation}
where %$C$ depends on the apriori assumptions and on  
%$M=\text{max}\left(\left\Vert \nabla_{g}a^{+}\right\Vert _{B_{1}^{+}},\left\Vert \nabla_{g}a^{-}\right\Vert _{B_{1}^{-}}\right)$, 
\[
T_{1}=\intop_{B_{2r}\setminus B_{r}}w^{2-2\tau}\left(r^{-2}\left|\nabla_{g}U\right|^{2}+r^{-4}U^{2}\right),
\]
and
\[
T_{2}=\intop_{B_{\frac{\overline{r}_{1}}{2}}\setminus B_{\frac{\overline{r}_{1}}{4}}^{\pm}}w^{2-2\tau}\left(\overline{r}_{1}^{-2}\left|\nabla_{g}U\right|^{2}+\overline{r}_{1}^{-4}U^{2}\right).
\]
Notice that in the above estimate, part of the integral of the term containing $\left|\nabla_{g}U\right|^{2}$ has been included in $T_1$ and $T_2$.

By \eqref{stimewl} and \eqref{sigma}, by using Caccioppoli inequality, we can write
\begin{equation}\label{T1}
    T_1\leq \frac{C}{ r^{2\tau+2}}\intop_{B_{2r}\setminus B_{r}^{\pm}}\left(r^{2}\left|\nabla_{g}U\right|^{2}+U^{2}\right)\leq \frac{C}{ r^{2\tau+2}}\intop_{B_{4r}}U^2
\end{equation}
and, in a similar way,
\begin{equation}\label{T2}
    T_2\leq \frac{C}{ \overline{r}_1^{2\tau+2}}\intop_{B_{\overline{r}_{1}}}U^{2}
\end{equation}

Now we estimate from below the first two terms in right hand side of \eqref{4.19} by using the definition of $\xi$, and we get
\begin{align}\label{rhs}
    \tau\intop_{B_1} w^{1-2\tau}\left|\nabla_{g}(\xi U)\right|^{2}+\tau^3\intop_{B_1}w^{-1-2\tau}(\xi U)^2\geq\nonumber\\
    \geq  \tau\intop_{B_{\frac{\overline{r}_1}{4}}\setminus B_{2r}} w^{1-2\tau}\left|\nabla_{g} U\right|^{2}+\tau^3\intop_{B_{\frac{\overline{r}_1}{4}}\setminus B_{2r}}w^{-1-2\tau} U^2
\end{align}

From \eqref{4.19}, \eqref{4.21}, \eqref{T1}, \eqref{T2} and \eqref{rhs} we have

\begin{align}\label{4.25}
    C\left(\intop_{B_{\frac{\overline{r}_{1}}{4}}\setminus B_{2r}}w^{2-2\tau}\left|\nabla_{g}U\right|^{2} +\frac{1}{ r^{2\tau+2}}\intop_{B_{4r}}U^2+\frac{1}{\overline{r}_1^{2\tau+2}}\intop_{B_{\overline{r}_{1}}}U^{2}\right)\geq\nonumber\\
    \geq \tau\intop_{B_{\frac{\overline{r}_1}{4}}\setminus B_{2r}} w^{1-2\tau}\left|\nabla_{g} U\right|^{2}+\tau^3\intop_{B_{\frac{\overline{r}_1}{4}}\setminus B_{2r}}w^{-1-2\tau} U^2+\tau^2r\intop_{B_1}w^{-2\tau-3}(\xi U)^2
\end{align}
Notice that, for $\tau$ large enough ($\tau>\max\{\tau_1, C\frac{\overline{r}_1}{4}\}$) we have
\begin{equation*}
 \tau\intop_{B_{\frac{\overline{r}_1}{4}}\setminus B_{2r}} w^{1-2\tau}\left|\nabla_{g} U\right|^{2}\geq C\intop_{B_{\frac{\overline{r}_{1}}{4}}\setminus B_{2r}}w^{2-2\tau}\left|\nabla_{g}U\right|^{2}   
\end{equation*}
so that \eqref{4.25} becomes
\begin{align}\label{4.26}
%\frac{\tau}{2}\intop_{B_{\frac{\overline{r}_1}{4}}\setminus B_{2r}} w^{1-2\tau}\left|\nabla_{g} U\right|^{2}+
\tau^3\intop_{B_{\frac{\overline{r}_1}{4}}\setminus B_{2r}}w^{-1-2\tau} U^2+\tau^2r\intop_{B_1}w^{-2\tau-3}(\xi U)^2\leq\nonumber\\
    \leq C\left( \frac{1}{ r^{2\tau+2}}\intop_{B_{4r}}U^2+\frac{1}{\overline{r}_1^{2\tau+2}}\intop_{B_{\overline{r}_{1}}}U^{2}\right)
\end{align}
%{\color{magenta} Siccome il pezzo con il graidente poi lo buttiamo via, non si può prendere direttamente $\tau$ in modo da poterlo buttare via, vale a dire  prendere 
%\begin{equation*}
 %\tau\intop_{B_{\frac{\overline{r}_1}{4}}\setminus B_{2r}} w^{1-2\tau}\left|\nabla_{g} U\right|^{2}\geq C\intop_{B_{\frac{\overline{r}_{1}}{4}}\setminus B_{2r}}w^{2-2\tau}\left|\nabla_{g}U\right|^{2} %  
%\end{equation*}
%}

Now, since $0<4r<R\leq\frac{\overline{r}_{1}}{4}$ and $\tau\geq 1$ and by the definition of $\xi$, we have 
\begin{equation}\label{4.27}
\tau^2r\intop_{B_1}w^{-2\tau-3}(\xi U)^2\geq \frac{1}{ (4r)^{2\tau+2}} \intop_{B_{4r}\setminus B_{2r}}U^2
\end{equation}
and
\begin{equation}\label{4.28}
\tau^3\intop_{B_{\frac{\overline{r}_1}{4}}\setminus B_{2r}}w^{-1-2\tau} U^2\geq \frac{1}{R^{2\tau+1}}\intop_{B_{R}\setminus B_{2r}} U^{2}.
\end{equation}

Combining \eqref{4.26}, \eqref{4.27}, and \eqref{4.28} 
 we have
\begin{align}\label{4.29}
\frac{1}{R^{2\tau+1}}\intop_{B_{R}\setminus B_{2r}} U^{2}& +\frac{1}{(4r)^{2\tau+2}} \intop_{B_{4r}\setminus B_{2r}}U^2\nonumber \\
\leq & C\left(\frac{1}{r^{2\tau+2}}\intop_{B_{4r}}U^2+\frac{1}{\overline{r}_1^{2\tau+2}}\intop_{B_{\overline{r}_{1}}}U^{2}\right)
\end{align}

Now, by adding the quantity $\left(R^{-1-2\tau}+
R\left(4r\right)^{-2\tau-2}\right)\intop_{B_{2r}}U^{2}$ 
to both sides of \eqref{4.29} we finally have  \eqref{tesilemma}
\end{proof}

\begin{theorem}\label{Th-Doubling}
(Pieciewise doubling inequality) Let $U$ be the solution of \eqref{eq:piecewiseeq_pb}.

There exists some constant $C$ that depends only on
$\lambda$, $\Lambda$, $\gamma_0$ and $\|\tilde{a}^\pm\|_{C^{0,1}(B_1^\pm)}$ such that, for every $0<r<\frac{\overline{r}_{1}}{16}$, 
\begin{equation}\label{doubfin}
\intop_{B_{4r}}U^{2}\leq \frac{CN^{3}}{\overline{r}_{1}^{3}}\intop_{B_{2r}}U^{2},
\end{equation}

where $N=\frac{\intop_{B_{\overline{r}_{1}}}U^{2}}{\intop_{B_{\overline{r}_{1}/4}}U^{2}}$
represent the largest frequency in the half ball.
\end{theorem}
\begin{proof}
Taking $\tau>\tau_{2}$ and take $R=\frac{\overline{r}_{1}}{4}$ in \eqref{tesilemma},
we have 
\begin{align}
\left(\frac{4}{\overline{r}_{1}}\right)^{1+2\tau}\intop_{B_{\overline{r}_{1}/4}}U^{2}+\frac{\overline{r}_{1}}{4}\frac{1}{\left(4r\right)^{2\tau+2}}\intop_{B_{4r}}U^{2}\nonumber \\
\leq C\left(\frac{1}{r^{2\tau+2}}\intop_{B_{2r}^{\mathcal{I}}}U^{2}+\frac{1}{\overline{r}_{1}^{2\tau+2}}\intop_{B_{\overline{r}_{1}}}U^{2}\right).\label{eq:eq1indoubling}
\end{align}

If we take $\tau=\overline{\tau}=\max\{\tau_2, \frac{1}{2}\log_4\frac{CN}{\overline{r}_1}\}$ we have 

\begin{equation*}
\left(\frac{4}{\overline{r}_{1}}\right)^{1+2\overline{\tau}}\intop_{B_{\overline{r}_{1}/4}}U^{2}\geq C\frac{1}{\overline{r}_{1}^{2\overline{\tau}+2}}\intop_{B_{\overline{r}_{1}}}U^{2},
\end{equation*}
so that \eqref{eq:eq1indoubling} becomes
\begin{align}
\intop_{B_{4r}}U^{2}
\leq \frac{C4^{2\overline{\tau}+3}}{\overline{r}_1}\intop_{B_{2r}^{\mathcal{I}}}U^{2}.\label{eq:eq2indoubling}
\end{align}
that is \eqref{doubfin}.
\end{proof}

\bibliographystyle{plain}

\begin{thebibliography}{99}
	\bibitem{A-E} V. Adolfsson L. Escauriaza, $C^{1,\alpha}$ domains and unique continuation at the boundary, Comm. Pure Appl. Math. (1997), 935--969. 
	
	%\bibitem{MRV2007} A. Morassi, E. Rosset and S. Vessella, \textit{Size estimates for inclusions in an elastic plate by boundary measurements}, Indiana U. Math. J. 56 (5) (2007), 2325--2384.
	
	
	%\bibitem{AMR2002} G. Alessandrini, A. Morassi, E. Rosset. Detecting an inclusion in an elastic body by boundary measurements. SIAM J. Math. Anal. 33 (2002) 1247--1268.
	
	\bibitem{AMR2003} G. Alessandrini, A. Morassi, E. Rosset. Size estimates, in ``Inverse problems: theory and applications'', Contemp. Math.  333 (2003) 1--33.	
	
%% non cit	
%	\bibitem{A-R-R-V} G. Alessandrini, L. Rondi, E. Rosset, S. Vessella. The stability for the Cauchy problem for elliptic equations. Inverse Problems 25 (2009) 1--47.
%%
	
	%\bibitem{AR1998} G. Alessandrini, E. Rosset. The inverse conductivity problem with one measurement: bounds on the size of the unknown object. SIAM J. Appl. Math. 58 (1998) 1060--1071.
	
	\bibitem{ARS2000} G. Alessandrini, E. Rosset, J.K. Seo. Optimal size estimates for the inverse conductivity problem with one measurement. Proc. Amer. Math. Soc. 128 (2000) 53--64. 
	
	\bibitem{Car1} T. Carleman, Sur les syst\`{e}mes lin\'{e}aires aux d\`{e}riv\`{e}es
	partielles du primier ordre \`{a} deux variables, C. R. Acad. Sci.
	Paris, 19, (1933), 471--474.
	
	\bibitem{Car2} T. Carleman, Sur un probl\`{e}me d'unicit\'{e} pour
	les syst\`{e}mes d'\`{e}quations aux d\`{e}riv\`{e}es partielles
	\`{a} deux variables ind\'{e}pendentes, Ark. Mat. Astr. Fys., 26B
	(1939), 1--9.
	
	\bibitem{Coi-Fe} R.R. Coifman, C.L. Fefferman. Weighted norm inequalities for maximal functions and singular integrals. Stud. Math. 51 (1974) 241--250. 
	
	\bibitem{DcFLVW} M. Di Cristo, E. Francini, C.--L. Lin, S. Vessella, J.--N- Wang, Carleman estimate for second order elliptic equations with Lipschitz leading coefficients and jump at an interface, J. Math. Pures Appl. (9) 108 (2), (2017), 163--206.
	
	\bibitem{E-V} L. Escauriaza and S. Vessella,  Optimal Three Cylinder Inequalities for Solutions to Parabolic Equations with Lipschitz Leading Coefficients (Inverse Problems Theory and Applications Contemporary Mathematics vol 333) ed G. Alessandrini and G. Uhlmann (Providence, RI: American  Mathematical Society), (2003), 79--87.
	
	\bibitem{FrLVW} E. Francini, C.-L. Lin, S. Vessella, J.-N. Wang, Three--region inequalities for the second order elliptic equation with discontinuous coefficients and size estimate, J. Differ. Equ. 261 (10), (2016), 5306--5323.
	
	\bibitem{Fr-Ve-Wa} E. Francini, S. Vessella, J.--N. Wang, Carleman estimate for complex second order elliptic operators with discontinuous Lipschitz coefficients, J. Spectr. Theory 12 (2), (2022), 535--571.
	
	\bibitem{Ga-Li} N. Garofalo, F. Lin. Monotonicity properties of variational integrals, $A_p$ weights and unique continuation. Indiana Univ. Math. J. 35 (1986) 245--268.
	
	\bibitem{G-Tr} D. Gilbarg, N.S. Trudinger. Elliptic Partial Differential Equations of Second Order. Springer, New York, 1983.
	
	\bibitem{HO63}  L. H\"{o}rmander. Linear Partial Differential Operators. Springer, New York, 1963.
	
	\bibitem{I98} V. Isakov. Inverse Problems for Partial Differential Equations. Applied Mathematical Sciences, vol. 127. Springer-Verlag, New York, 1998.
	
	
	\bibitem{La-Ur} O. Ladyzhensaya and N. Ural'teva, Linear and quasi--linear elliptic equations, Academic Press, New York, 1968.
	
	\bibitem{Le-Lern} J. Le Rousseau, N. Lerner, Carleman estimates for anisotropic elliptic operators with jumps at an interface, Anal. PDE 6 (7), (2013), 1601--1648.
	
	\bibitem{Ler} N. Lerner. Carleman inequalities. An introduction and more. Grundlehren der mathematischen Wissenschaften, Springer Verlag, Berlin, 2019.
	
%%non cit
%	\bibitem{LiNW} C.-L. Lin, S. Nagayasu, J.-N. Wang. Quantitative uniqueness for the power of the Laplacian with singular coefficients. Ann. Sc. Norm. Super. Pisa Cl. Sci. (5) 10 (2011) 513--529.
%%	
	\bibitem{Mand} N. Mandache N, On a counterexample concerning unique continuation
	for elliptic equations in divergence form,  Mat. Fiz. Anal. Geom. 3, (1996), 308-31.
	
	\bibitem{Mil} K. Miller, Nonunique continuation for uniformly parabolic and elliptic equations in selfadjoint divergence form with Hölder continuous coefficients. Arch. Rational Mech. Anal. 54, (1974), 105–117.
	
	
%% non cit	
%	\bibitem{MRV2007} A. Morassi, E. Rosset, S. Vessella. Size estimates for inclusions in an elastic plate by boundary measurements. Indiana Univ. Math. J. 56 (2007) 2325--2384.
%%	
	%\bibitem{MRV2009} A. Morassi, E. Rosset, S. Vessella. Detecting general inclusions in elastic plates. Inverse Problems 25 (2009) 1--14.
	
	
	%\bibitem{MRV2018} A. Morassi, E. Rosset, S. Vessella. Size estimates for fat inclusions in an isotropic Reissner-Mindlin plate. Inverse Problems 34 (2018) Paper 025001. 
	
	%\bibitem{mrv:mat-cat} A. Morassi, E. Rosset, S. Vessella. Doubling inequality at
	%the boundary for the Kirchhoff-Love plate's equation with Dirichlet
	%conditions. Le Matematiche 75 (2020) 27--55.
	
	
	
	%\bibitem{Pl} A. Pl\u{\i}s, On non-uniqueness in Cauchy problem for an elliptic
	%second order differential equation, Bull. Acad. Pol. Sci. S\'{e}r.
	%Sci. Math. Astron. Phys. 11, (1963), 95 -- 100.
	
	
	\bibitem{Ves} S. Vessella, Unique Continuation Properties for Partial Differential Equations--Introduction to the Stability Estimates for Inverse Problems, Birkh\"{a}user, Springe Nature Gewerbestrasse 11, 6330 Cham, Switzerland.  
	
		
\end{thebibliography}

\end{document}